\input amstex
\magnification=1200
\loadmsbm
\loadeufm
\loadeusm
\UseAMSsymbols
\input amssym.def

\font\BIGtitle=cmr10 scaled\magstep3
\font\bigtitle=cmr10 scaled\magstep1
\font\boldsectionfont=cmb10 scaled\magstep1
\font\section=cmsy10 scaled\magstep1

\def\scr#1{{\fam\eusmfam\relax#1}}
\def\scrA{{\scr A}}
\def\scrB{{\scr B}}
\def\scrC{{\scr C}}
\def\scrD{{\scr D}}
\def\scrE{{\scr E}}
\def\scrF{{\scr F}}
\def\scrG{{\scr G}}

\def\scrH{{\scr H}}

\def\scrL{{\scr L}}

\def\scrM{{\scr M}}
\def\scrN{{\scr N}}
\def\scrO{{\scr O}}

\def\scrS{{\scr S}}
\def\scrU{{\scr U}}

\def\scrT{{\scr T}}
\def\scrV{{\scr V}}
\def\scrX{{\scr X}}
\def\scrY{{\scr Y}}
\def\scrW{{\scr W}}

\def\gr#1{{\fam\eufmfam\relax#1}}

\def\grC{{\gr C}}	
\def\grD{{\gr D}}	
\def\grE{{\gr E}}	
	
\def\grG{{\gr G}}	
\def\grH{{\gr H}}

\def\grK{{\gr K}}	
\def\grL{{\gr L}}

\def\grO{{\gr O}}

	\def\grs{{\gr s}}
	\def\grt{{\gr t}}

\def\db#1{{\fam\msbfam\relax#1}}

 \def\dbB{{\db B}}
 \def\dbD{{\db D}}
\def\dbE{{\db E}} \def\dbF{{\db F}}
\def\dbG{{\db G}}

 \def\dbN{{\db N}}
 \def\dbP{{\db P}}
\def\dbQ{{\db Q}} 
 \def\dbT{{\db T}}

 \def\dbZ{{\db Z}}

\def\eps{{\varepsilon}}

\def\im{\text{Im}}

\def\Ker{\text{Ker}}
\def\der{\text{der}}

\def\GL{\text{GL}}
\def\Hom{\text{Hom}}
\def\End{\text{End}}
\def\Spec{\text{Spec}}
\def\Spf{\text{Spf}}

\def\Tr{\text{Tr}}

\def\Lie{\text{Lie}}

\def\leaderfill{\leaders\hbox to 1em
     {\hss.\hss}\hfill}
\def\nspace{\lineskip=1pt\baselineskip=12pt\lineskiplimit=0pt}

\def\finishproclaim{\par\rm
     \ifdim\lastskip<\medskipamount\removelastskip
     \penalty55\medskip\fi}
\def\endproof{$\hfill \square$}
\def\proof{\par\noindent {\it Proof:}\enspace}
\def\references#1{\par
  \centerline{\boldsectionfont References}\medskip
     \parindent=#1pt\nspace}
\def\Ref[#1]{\par\hang\indent\llap{\hbox to\parindent
     {[#1]\hfil\enspace}}\ignorespaces}
\def\Item#1{\par\smallskip\hang\indent\llap{\hbox to\parindent
     {#1\hfill$\,\,$}}\ignorespaces}
\def\ItemItem#1{\par\indent\hangindent2\parindent
     \hbox to \parindent{#1\hfill\enspace}\ignorespaces}

\def\arrowsim{\,\smash{\mathop{\to}\limits^{\lower1.5pt
  \hbox{$\scriptstyle\sim$}}}\,}

\def\doublemaprights#1#2#3#4{\raise3pt\hbox{$\mathop{\,\,\hbox to     
#1pt{\rightarrowfill}\kern-30pt\lower3.95pt\hbox to
     #2pt{\rightarrowfill}\,\,}\limits_{#3}^{#4}$}}

\def\rightcapdownarrow{\raise9pt\hbox{$\ssize\cap$}\kern-7.75pt
     \Big\downarrow}

\def\rcapmapdown#1{\rightcapdownarrow\kern-1.0pt\vcenter{
     \hbox{$\scriptstyle#1$}}}

\def\rmapdown#1{\Big\downarrow\kern-1.0pt\vcenter{
     \hbox{$\scriptstyle#1$}}}
\def\rightsubsetarrow#1{{\ssize\subset}\kern-4.5pt\lower2.85pt
     \hbox to #1pt{\rightarrowfill}}
\def\longtwoheadedrightarrow#1{\raise2.2pt\hbox to #1pt{\hrulefill}
     \!\!\!\twoheadrightarrow}

\def\Hom{\operatorname{\hbox{Hom}}}

\NoBlackBoxes
\parindent=25pt
\document
\footline={\hfil}
\footline={\hss\tenrm \folio\hss}
\pageno=1
\bigskip

\bigskip\bigskip\bigskip\bigskip
\noindent 
\centerline{\BIGtitle Deformation subspaces of $p$-divisible groups}
\bigskip\bigskip
\centerline{\BIGtitle  as formal Lie groups associated to $p$-divisible groups}
\bigskip\bigskip
\centerline{\bigtitle Adrian Vasiu, Binghamton University, U.S.A.}
\bigskip
\centerline{April 20, 2010, accepted in final form for publ. in J. Alg. Geom.}

\bigskip\noindent
{\bf ABSTRACT.} Let $k$ be an algebraically closed field of characteristic $p>0$. Let $D$ be a $p$-divisible group over $k$ which is not isoclinic. Let $\scrD$ (resp. $\scrD_k$) be the formal deformation space of $D$ over $\Spf(W(k))$ (resp. over $\Spf(k)$). We use axioms to construct formal subschemes $\scrG_k$ of $\scrD_k$ that: (i) have canonical structures of formal Lie groups over $\Spf(k)$ associated to $p$-divisible groups over $k$, and (ii) give birth, via all geometric points $\Spf(K)\to\scrG_k$, to $p$-divisible groups over $K$ that are isomorphic to $D_K$. We also identify when there exist formal subschemes $\scrG$ of $\scrD$ which lift $\scrG_k$ and which have natural structures of formal Lie groups over $\Spf(W(k))$ associated to $p$-divisible groups over $W(k)$. Applications to Traverso (ultimate) stratifications are included as well. 

\bigskip\noindent
{\bf KEY WORDS:} $p$-divisible groups, formal deformation spaces, formal Lie groups, $F$-crystals, Lie algebras, affine group schemes, and stratifications. 

\bigskip\noindent
{\bf MSC 2000:} 11G10, 11G18, 14F30, 14G35, 14L05, and 20G25.

\bigskip\bigskip
\noindent
{\boldsectionfont 1. Introduction}
\bigskip

Let $p\in\dbN$ be a prime. Let $k$ be an algebraically closed field of characteristic $p$. Let $W(k)$ be the ring of Witt vectors
 with coefficients in $k$. Let $B(k)$ be the field of fractions of $W(k)$. Let $\sigma:=\sigma_k$ be
the Frobenius automorphism of $W(k)$ and $B(k)$ induced from $k$. Let $c,d\in\dbN\cup\{0\}$ be such that $r:=c+d>0$. Let $D$ be a {\it $p$-divisible group} over $k$ of codimension $c$ and dimension $d$; its height is $r$. Let $\dbD(D)=(M,\phi,\vartheta)$ be the {\it crystalline (contravariant) Dieudonn\'e module} of $D$ in the sense of [BBM, Def. 3.1.5], where $\phi$ and $\vartheta$ are the $F$
and $V$ of loc. cit. We recall that $M$ is a free $W(k)$-module of rank $r$ and that $\phi:M\to M$ is a $\sigma$-linear endomorphism such that we have $pM\subseteq \phi(M)$. The Verschiebung map $\vartheta:M\to M$ of $\phi$ is a $\sigma$-linear endomorphism such that we have $\vartheta\circ\phi=\phi\circ\vartheta=p1_M$. In what follows we will simply call $(M,\phi)$ the Dieudonn\'e module of $D$. We denote also by $\phi$ the
$\sigma$-linear automorphism of $\End(M[{1\over p}])$ that takes
$e\in \End(M[{1\over p}])$ to $\phi(e):=\phi\circ
e\circ\phi^{-1}=\vartheta^{-1}\circ e\circ\vartheta$. 

If $\Spec(A)$ is a local, complete, noetherian scheme, let $\Spf(A)$ be the affine formal scheme defined by the maximal ideal of $A$. Let $R:=W(k)[[x_1,\ldots,x_{cd}]]$. The {\it formal deformation space} $\scrD$ of $D$ over $\Spf(W(k))$ is representable by $\Spf(R)$, cf. [Il, Cor. 4.8 (i)]. Throughout the paper we fix an identification $\scrD=\Spf(R)$. Let $\scrD_k:=\Spf(R/pR)$; we identify $\scrD_k$ with the formal deformation space of $D$ over $\Spf(k)$. The categories of $p$-divisible groups over $\Spec(R)$ and $\Spf(R)$ are canonically identified, cf. [Me, Ch. II, Lemma 4.16]. Therefore there exists a unique $p$-divisible group $\grD$ over $\Spec(R)$ that gives birth to the universal $p$-divisible group $\grD^{\text{f}}$ over $\scrD=\Spf(R)$. Both $\grD^{\text{f}}$ and $\grD$ are versal deformations of $D$, cf. [Il, Cor. 4.8 (ii)]. 

A {\it formal Lie group} $\scrF$ (resp. $\scrF_k)$ over $\Spf(W(k))$ (resp. over $\Spf(k)$) is a contravariant functor from the category of local, artinian $W(k)$-schemes (resp. $k$-schemes) of residue field $k$ into the category of groups which, when viewed as a functor into the category of sets, is representable by $\Spf({S})$ (resp. by $\Spf({S}/pS)$), where $S=W(k)[[x_1,\ldots,x_n]]$ for some number $n\in\dbN\cup\{0\}$ called the dimension of $\scrF$ (resp. of $\scrF_k)$ (see [Me, Ch. 2, Def. 1.1.5] and [Fo]). We recall that the category of connected $p$-divisible groups over $W(k)$ (resp. over $k$) is equivalent to the category of divisible commutative formal Lie groups over $\Spf(W(k))$ (resp. over $\Spf(k)$), cf. [Ta, Sect. 2, (2.2), Prop. 1]. This equivalence of categories is achieved via the association on objects $\vartriangle\to Hom((\dbZ_p)_{W(k)},\vartriangle)$ (resp. $\vartriangle\to Hom((\dbZ_p)_k,\vartriangle)$). The short exact sequence of group schemes $0\to (\dbZ_p)_{W(k)}\to (\dbQ_p)_{W(k)}\to (\dbQ_p/\dbZ_p)_{W(k)}\to 0$ over $W(k)$ (resp. its analogue over $k$) gives birth to a coboundary functor 
$$Hom((\dbZ_p)_{W(k)},\vartriangle)\to Ext^1((\dbQ_p/\dbZ_p)_{W(k)},\vartriangle)$$ 
(resp. $Hom((\dbZ_p)_k,\vartriangle)\to Ext^1((\dbQ_p/\dbZ_p)_k,\vartriangle)$) which is an isomorphism of formal Lie groups over $\Spf(W(k))$ (resp. over $\Spf(k)$) (see [Va4, Lemma 5.3.1.2] for this well known fact) and which is functorial in $\vartriangle$ (this is a standard cohomological property). 

The classical {\it Serre--Tate ordinary theory} shows that if $D$ is ordinary, then $\scrD$ has a canonical structure of a formal torus over $\Spf(W(k))$ (see [De, Subsect. 2.1 A)] and [Ka2, Thm. A3.1]). This theory was generalized in [Va4] to abstract and geometric contexts pertaining to integral canonical models of Shimura varieties of Hodge type. Roughly, in [Va4, Chs. 5 and 9] it is shown that in many cases (some of them are described in Example 1.8), $\scrD$ has formal subschemes $\scrG$ that have the following two main properties:

\medskip
{\bf (i)} they have canonical structures of formal Lie groups over $\Spf(W(k))$ which are isomorphic to formal Lie groups over $\Spf(W(k))$ associated to $p$-divisible groups over $W(k)$ whose special fibres have Dieudonn\'e modules that {\it sit} naturally inside the pair $(\End(M),\phi)$;

\smallskip
{\bf (ii)} the restriction of $\grD^{\text{f}}$ to $\scrG\cap\scrD_k$ is geometrically constant in the following sense. 

\medskip\smallskip\noindent
{\bf 1.1. Definition.} Let $\grC$ be a $p$-divisible group over $\Spec(k[[x_1,\ldots,x_n]])$ that lifts $D$. Let $\grC^{\text{f}}$ be the $p$-divisible group over $\Spf(k[[x_1,\ldots,x_n]])$ defined by $\grC$. We say that either $\grC$ or $\grC^{\text{f}}$ is {\it geometrically constant}, if the following property holds:
For each algebraically closed field $K$ that contains $k$ and for every morphism $f:\Spec(K)\to\Spec(k[[x_1,\ldots,x_n]])$, the $p$-divisible group $f^*(\grC)$ is isomorphic to $D_K$.

\medskip
The goal of the paper is to put forward {\it basic principles} that govern the phenomena described by properties (i) and (ii) above. We will use {\it axioms} to identify: (a) formal subschemes $\scrG_k$ of $\scrD_k$ for which $\scrG_k\times_{\scrD_k} \grD^{\text{f}}$ is geometrically constant and which have canonical structures of formal Lie groups over $\Spf(k)$ isomorphic to formal Lie groups over $\Spf(k)$ associated to $p$-divisible groups over $k$, and (b) when there exist formal subschemes $\scrG$ of $\scrD$ which lift $\scrG_k$ and which have natural structures of formal Lie groups over $\Spf(W(k))$ isomorphic to formal Lie groups over $\Spf(W(k))$ associated to $p$-divisible groups over $W(k)$. 

The importance of the smooth formal subschemes $\scrG_k$ and $\scrG$ stems from the study of smooth integral models of Shimura varieties of Hodge type over local, \'etale $\dbZ_{(p)}$-algebras. Each such smooth integral model $\scrM$ is a moduli space of polarized abelian schemes which are endowed with symplectic similitude level structures and with families of specializations of Hodge cycles (see [Mi1, Sect. 2], [Mi2, Sects. 3 and 4], [Va1, Sect. 4], etc.). The special fibre $\scrS$ of $\scrM$ is endowed with many types of {\it stratifications} as well as with certain {\it groups of Hecke operators} which leave invariant the corresponding strata. A very important stratification of $\scrS$ is the {\it Traverso (ultimate) stratification} which can be defined similarly to either [Va2, Thm. 5.3.1 and Subsubsect. 5.3.2] or [Va3, Subsubsect. 4.2.4] and whose strata parametrize certain isomorphism classes of Dieudonn\'e modules endowed with extra structures (i.e., with a group as in [Va2, Def. 1.1]). In practice, the $\scrG_k$'s are {\it slices} of the {\it local Traverso strata} i.e., are slices of the smooth formal subschemes of $\scrD_k$ that are defined by strata of the most general Traverso stratifications and that are invariant under large subgroups of Hecke operators (see Subsection 5.2 for local Traverso strata defined abstractly). Moreover, each $\scrG$ is a good lift of $\scrG_k$ over $\Spf(W(k))$ which is supposed to be also invariant under large subgroups of Hecke operators and which (for instance, see Example 1.8) can be often used to identify canonical or quasi-canonical lifts of $D$ (endowed with extra structures) to finite, discrete valuation ring extensions of $W(k)$. 

Each {\it new structure} on $\scrD_k$ (like the ones provided by the formal Lie group structures on $\scrG_k$'s) represents a {\it new tool} towards the study of groups of Hecke operators on $\scrS$. Such tools play key roles in identifying Shimura subvarieties of Shimura varieties of Hodge type (see [Mo]), in computing different cohomology groups of simple Shimura varieties (see [HT]), in studying the Zariski closures of different Hecke orbits (see [Ch1]), etc.  

In Subsection 1.2 we introduce the language of Newton polygon slopes. In Subsection 1.3 we recall standard descriptions of the tangent space of $\scrD_k$. In Subsection 1.4 we present four axioms. Our main results are stated in  Subsections 1.5 and 1.6. Subsections 1.7 to 1.9 contain examples. Subsection 1.10 outlines very briefly the main problems left unsolved. Subsection 1.11 details on the contents of the paper.

\medskip\smallskip\noindent
{\bf 1.2. Slopes.} Dieudonn\'e's classification of $F$-isocrystals over $k$ (see [Di, Thms. 1 and 2], [Ma, Ch. 2, Sect. 4], etc.) implies that we have a direct sum decomposition 
$$M[{1\over p}]=\oplus_{\alpha\in\scrW} W(\alpha)$$
that is stable under $\phi$ and that has the property that all Newton polygon slopes of $(W(\alpha),\phi)$ are $\alpha$; here $\scrW$ is the set of Newton polygons slopes of $(M,\phi)$. We recall that $\scrW\subseteq\dbQ\cap [0,1]$. We also recall that if $m_{\alpha}\in\dbN$ is the smallest number such that $m_{\alpha}\alpha\in\dbN\cup\{0\}$, then there exists a $B(k)$-basis for $W(\alpha)$ which is formed by elements fixed by $p^{-m_{\alpha}\alpha}\phi^{m_{\alpha}}$. One says that $D$ is {\it isoclinic} if the set $\scrW$ has only one element. We consider the direct sum decomposition $\End(M[{1\over p}])=L_+\oplus L_0\oplus L_-$ 
that is stable under $\phi$ and such that all Newton polygon slopes of $(L_+,\phi)$ are positive, all Newton polygon slopes of $(L_0,\phi)$ are $0$, and all Newton polygon slopes of $(L_-,\phi)$ are negative. Let 
$$V_+:=\End(M)\cap L_+\;\;\text{and}\;\; V_-:=\End(M)\cap L_-.$$ 
\indent
We have three direct sum decompositions $L_+=\oplus_{\alpha,\beta\in\scrW,\,\alpha<\beta} \Hom(W(\alpha),W(\beta))$, $L_0=\oplus_{\alpha\in\scrW} \End(W(\alpha))$, and $L_-=\oplus_{\alpha,\beta\in\scrW,\,\alpha<\beta} \Hom(W(\beta),W(\alpha))$. 

\medskip\smallskip\noindent
{\bf 1.3. Tangent space of $\scrD_k$.} Let $D^{\text{t}}$ be the Cartier dual of $D$. Let $t_D$ and $t_{D^{\text{t}}}$ be the tangent spaces of $D$ and $D^{\text{t}}$ (respectively) at the identity element. Let $\bar M:=M/pM$. Let $\bar\phi:\bar M\to\bar M$ be the reduction modulo $p$ of $\phi$. 
Let $\bar F^1:=\Ker(\bar\phi)$; it is a $k$-vector subspace of $\bar M$ that has dimension $d$. Let $F^0(\End(\bar M)):=\{e\in\End(\bar M)|e(\bar F^1)\subseteq\bar F^1\}$. Let 
$$\scrT:=\End(\bar M)/F^0(\End(\bar M)).$$
If $\bar x\in\Hom(\bar F^1,\bar M/\bar F^1)$ and if $\bar Y\in\End(\bar M)$ is an element whose restriction to $\bar F^1$ defines a $k$-linear map $\bar X:\bar F^1\to\bar M$ that lifts $\bar x$, then the image $\bar y$ of $\bar Y$ in $\scrT$ depends only on $\bar x$. The association $\bar x\to\bar y$ is $k$-linear and defines a canonical identification $\Hom(\bar F^1,\bar M/\bar F^1)=\scrT$.  

Identifying $\bar F^1$ with the dual $k$-vector space of $t_D$, the tangent space $t_{D^{\text{t}}}$ gets identified with $\bar M/\bar F^1$. Thus we can identify $t_{D^{\text{t}}}\otimes_k t_D=\bar M/\bar F^1\otimes_k \Hom(\bar F^1,k)=\Hom(\bar F^1,\bar M/\bar F^1)$. As the tangent space of $\scrD_k$ is canonically identified with $t_{D^{\text{t}}}\otimes t_D$ (cf. [Il, Cor. 4.8 (iii)]), it is also canonically identified with $\Hom(\bar F^1,\bar M/\bar F^1)$ and therefore also with $\scrT$.  

\medskip\noindent
{\bf 1.3.1. The map $\nu$, lifts, and cocharacters.} Let 
$$\nu:\End(M)\twoheadrightarrow \scrT$$ 
be the natural $W(k)$-epimorphism. Let $F^1$ be a direct summand of $M$ that lifts $\bar F^1$; its rank is $d$. The triple $(M,F^1,\phi)$ is called a {\it filtered Dieudonn\'e module} over $k$. One also calls $F^1$ or $(M,F^1,\phi)$ as a {\it lift} of $(M,\phi)$. Let $\mu:\dbG_m\to \pmb{\GL}_M$ be the inverse of the {\it canonical split cocharacter} of $(M,F^1,\phi)$ defined in [Wi, p. 512]. We have a unique direct sum decomposition $M=F^1\oplus F^0$ such that $\dbG_m$ acts through $\mu$ trivially on $F^0$ and via the inverse of the identity character of $\dbG_m$ on $F^1$. We have a direct sum decomposition
$$\End(M)=\Hom(F^0,F^1)\oplus\End(F^1)\oplus\End(F^0)\oplus \Hom(F^1,F^0)$$ 
of $W(k)$-modules. Let $\mu_L:\dbG_m\to \pmb{\GL}_{\End(M)}$ be the cocharacter defined naturally by $\mu$; the group $\dbG_m$ acts via the inverse of the identity character of $\dbG_m$ on $\Hom(F^0,F^1)$, acts trivially on $\End(F^1)\oplus\End(F^0)$, and acts via the identity character of $\dbG_m$ on $\Hom(F^1,F^0)$. We have $\phi({1\over p}F^1+M)=M$. Thus we can write $\phi=\sigma_{\phi}\mu({1\over p})$, where $\sigma_{\phi}:M\arrowsim M$ is a $\sigma$-linear automorphism of $M$. Therefore for $x\in\End(M)$, we have $\phi(x)\in\End(M)$ if and only if $x$ belongs to $\Hom(F^0,F^1)\oplus\End(F^1)\oplus\End(F^0)\oplus p\Hom(F^1,F^0)=\Ker(\nu)$. Thus:

\medskip
{\bf (*)}  For $x\in\End(M)$, we have $\phi(x)\notin \End(M)$ if and only if $\nu(x)\neq 0$. 

\medskip
If $E_{+0}$ is a $W(k)$-submodule of $\End(M)$ such that the pair $(E_{+0},\phi)$ is a Dieudonn\'e module over $k$ (thus $E_{+0}\subseteq \End(M)\cap (L_+\oplus L_0)$), then we have $\nu(E_{+0})=0$ (cf. (*)). Therefore we can not use $E_{+0}$ to construct versal deformations of $D$. Thus, in order to construct versal deformations of $D$ via Dieudonn\'e modules, we will have to use $W(k)$-submodules $E_-$ of $V_-$ such that the pair $(E_-,p\phi)$ is a Dieudonn\'e module over $k$ and $\nu(E_-)\neq 0$. Accordingly, in the whole paper we will assume that $D$ is {\it not} isoclinic (i.e., we will assume that $V_+\neq 0$ and $V_-\neq 0$) and we will make use of the following list of axioms.

\medskip\smallskip\noindent
{\bf 1.4. Axioms.} For a $W(k)$-submodule $E$ of $V_-$, we consider the following two axioms:

\medskip
{\bf (i)} the pair $(E,p\phi)$ is the largest Dieudonn\'e module contained in $(E[{1\over p}]\cap V_-,p\phi)$; 

\smallskip
{\bf (ii)} inside the $W(k)$-algebra $\End(M)$, we have $E^2:=E\cdot E=0$.

\medskip
Axioms (i) and  (ii) imply that there exists an element $t\in L_0$ which is fixed by $\phi$ and such that we have $[E,t]=E$, cf. Lemma 3.2 (c.i). Thus axioms (i) and (ii) pertain to the possibility of using $E\oplus W(k)t$ to produce versal deformations of $D$ that give birth to the {\it universal extension} in characteristic $p$ of $(\pmb{\mu}_{p^{\infty}})_k$ by a $p$-divisible group over $k$ whose Dieudonn\'e module is $(E,p\phi)$; here we identify $\dbD((\pmb{\mu}_{p^{\infty}})_k)$ with $(W(k)t,p\phi,\phi^{-1})$. 

Let $D_{W(k)}$ be a $p$-divisible group over $W(k)$ that lifts $D$. Not to introduce extra notations, we assume that the direct summand $F^1$ of $M$ is the {\it Hodge filtration} of $D_{W(k)}$. Let $F^0(E):=E\cap (\End(F^1)\oplus \End(F^0))$ and $F^{-1}(E):=E\cap \Hom(F^1,F^0)$. For the pair $(E,D_{W(k)})$, we consider the following two additional axioms:

\medskip
{\bf (iii)} we have a direct sum decomposition $E=F^0(E)\oplus F^{-1}(E)$;

\smallskip
{\bf (iv)} we have a direct sum decomposition $E=\phi(F^0(E))\oplus p\phi(F^{-1}(E))$.

\medskip
Axiom (iii) implies that $\mu_L$ normalizes $E$. If the axioms (iii) and (iv) hold, then the triple $(E,F^0(E),p\phi)$ is a filtered Dieudonn\'e module over $k$. 

The two Basic Theorems below are proved in Subsections 3.3 and 3.5 (respectively).

\medskip\smallskip\noindent
{\bf 1.5. Basic Theorem A.} {\it Let $D$ be a $p$-divisible group over $k$ that is not isoclinic. Let $E$ be a $W(k)$-submodule of $V_-=\End(M)\cap L_-$ such that the axioms 1.4 (i) and (ii) hold for it. Let $\dbD^{-1}(E)$ be a $p$-divisible group over $k$ whose Dieudonn\'e module is isomorphic to $(E,p\phi)$. Then to $E$ corresponds canonically a formal subscheme $\scrG^E_k$ of $\scrD_k$ that has the following three properties:

\medskip
{\bf (i)} it has a canonical structure of a formal Lie group over $\Spf(k)$ isomorphic to the formal Lie group over $\Spf(k)$ associated to the Cartier dual $\dbD^{-1}(E)^{\text{t}}$ of $\dbD^{-1}(E)$ (thus the dimension of $\scrG^E_k$ is the codimension of $\dbD^{-1}(E)$);

\smallskip
{\bf (ii)} its tangent space is $\nu(E)$; 

\smallskip
{\bf (iii)} the restriction of $\grD^{\text{f}}$ to $\scrG^E_k$ is geometrically constant.}

\medskip
One thinks of $\scrG_k^E$ as keeping track of the {\it $E$-deformations} of $D$ (in characteristic $p$) in a sense that can be made very precise by using $F$-crystals (see Subsection 3.3).

\medskip\smallskip\noindent
{\bf 1.6. Basic Theorem B.} {\it Let $D$ and $E$ be as in the Basic Theorem A. Suppose that there exists a $p$-divisible group $D_{W(k)}$ over $W(k)$ that lifts $D$ and such that the axioms 1.4 (iii) to (iv) hold for the pair $(E,D_{W(k)})$. Then the following four properties hold:

\medskip
{\bf (i)} there exists a unique $p$-divisible group $\dbD^{-1}(E)_{W(k)}$ over $W(k)$ whose filtered Dieudonn\'e module is the triple $(E,F^0(E),p\phi)$;

\smallskip
{\bf (ii)} to $(E,D_{W(k)})$ corresponds naturally a formal subscheme $\scrG^E$ of $\scrD$ that lifts $\scrG^E_k$ and that has a natural structure of a formal Lie group over $\Spf(W(k))$ isomorphic to the formal Lie group over $\Spf(W(k))$ of the Cartier dual $\dbD^{-1}(E)_{W(k)}^{\text{t}}$ of $\dbD^{-1}(E)_{W(k)}$ (the formal Lie group structure induced on $\scrG^E_k$ is the canonical one of the Basic Theorem A);

\smallskip
{\bf (iii)} the origin of $\scrG^E$ is defined by $D_{W(k)}$ and the filtered Dieudonn\'e modules of $p$-divisible groups over $W(k)$ defined by $\Spf(W(k))$-valued points of $\scrG^E$, are of the form $(M,u(F^1),\phi)$, where $u\in 1_M+ pF^{-1}(E)$;

\smallskip
{\bf (iv)} if $D^\prime_{W(k)}$ is a $p$-divisible group over $W(k)$ defined by a morphism $\Spf(W(k))\to\scrG^E$ (so $D^\prime_{W(k)}$ lifts $D$), then the axioms 1.4 (iii) and (iv) hold for the pair $(E,D^\prime_{W(k)})$.}

\medskip\noindent
As (iv) and the notation $\scrG^E$ suggest, the main role of $D_{W(k)}$ is to fix the origin of $\scrG^E$.

Always there exists a short exact sequence $0\to D_2\to D\to D_1\to 0$ that depends on $E$ and such that $\scrG_k^E$ is naturally a formal Lie subgroup of the formal Lie group $Ext^1(D_1,D_2)$ over $\Spf(k)$ of extensions of $D_1$ by $D_2$ (cf. Example 3.4 and Corollary 3.7). But in general, the smallest Newton polygon slope of $D_2$ is not greater than the greatest Newton polygon slope of $D_1$ (cf. Example 1.9). 

If $D=D_1\times_k D_2$ is a direct sum of two isoclinic $p$-divisible groups of distinct Newton polygon slopes and if $E$ is the maximal $W(k)$-submodule of $V_-$ for which the axioms 1.4 (i) and (ii) hold, then an equivalent form of the Basic Theorem A and Example 3.4 were first proved in [Ch2] using complicated and lengthy computations on (triple) Cartier modules. Thus the main results of [Ch2] are encompassed in this paper by particular cases of Examples 1.7 and 3.4. There exists no particular analogue of Basic Theorem B in [Ch2]. Basic Theorems A and B are functorial in $E$ in the natural way, cf. Proposition 3.6. 

The following three examples illustrate key aspects of the Basic Theorems A and B. Example 4.3.2 illustrates in depth the functoriality of the  Basic Theorem A. 

\medskip\smallskip\noindent
{\bf 1.7. Example.} We consider a short exact sequence $0\to (M_1,\phi)\to (M,\phi)\to (M_2,\phi)\to 0$ of Dieudonn\'e modules over $k$ such that each Newton polygon slope of $(M_2,\phi)$ is greater than all Newton polygon slopes of $(M_1,\phi)$; to it corresponds naturally a short exact sequence $0\to D_2\to D\to D_1\to 0$ of $p$-divisible groups over $k$. Let $r_1$ and $r_2=r-r_1$ be the heights of $D_1$ and $D_2$ (respectively). The direct summand $\Hom(M_2,M_1)$ of $\End(M)$ is also a direct summand of $V_-$. Let $E$ be the largest $W(k)$-submodule of $\Hom(M_2,M_1)$ such that the pair $(E,p\phi)$ is the Dieudonn\'e module over $k$ of a $p$-divisible group $\dbD^{-1}(E)$. As $E^2=0$ and due to the choice of $(E,p\phi)$, the axioms 1.4 (i) and (ii) hold for $E$. 

Based on the Basic Theorem A we get that to $E$ corresponds a smooth formal subscheme $\scrG^E_k$ of $\scrD_k$ that has a canonical structure of a formal Lie group over $\Spf(k)$ isomorphic to the formal Lie group over $\Spf(k)$ associated to $\dbD^{-1}(E)^{\text{t}}$. Moreover $\scrG_k^E$ is a formal Lie subgroup of the formal Lie group $Ext^1(D_1,D_2)$ over $\Spf(k)$, cf. Corollary 3.7.

In this and the next paragraph we consider the particular case when $D_1$ and $D_2$ are isoclinic of Newton polygon slopes $\alpha_1$ and $\alpha_2$ (respectively). We have $\alpha_1<\alpha_2$, $\scrW=\{\alpha_1,\alpha_2\}$, and $V_-=\Hom(M_2,M_1)$. As $\Hom(M_2,M_1)=\Hom(M_2,W(k))\otimes_{W(k)} M_1$, the unique Newton polygon slope of $\dbD^{-1}(E)$ is $1-\alpha_2+\alpha_1$. Thus $\dbD^{-1}(E)^{\text{t}}$ is isoclinic of height $r_1r_2$, Newton polygon slope $\alpha_2-\alpha_1$, and dimension $r_1r_2(\alpha_2-\alpha_1)$. If $D$ is ordinary (i.e., if $(\alpha_1,\alpha_2)=(0,1)$), then $\dbD^{-1}(E)^{\text{t}}$ is isomorphic to $(\pmb{\mu}_{p^{\infty}})^{cd}_k$ and therefore the previous paragraph recovers the classical result of Serre--Tate. If $D$ is isomorphic to $D_1\times_k D_2$, then the previous paragraph recovers the main results of [Ch2]. 

To exemplify the Basic Theorem B, in this paragraph we will also assume that $D=D_1\times_k D_2$ (thus $M=M_1\oplus M_2$) and that there exist $W(k)$-bases $\{e_1,\ldots,e_{r_1}\}$ and $\{f_1,\ldots,f_{r_2}\}$ for $M_1$ and $M_2$ (respectively) and permutations $\pi_1\in S_{r_1}$ and $\pi_2\in S_{r_2}$ such that we have $\phi(e_i)=p^{a_i}e_{\pi_1(i)}$ and $\phi(f_j)=p^{b_j}f_{\pi_2(j)}$, where $i\in\{1,\ldots,r_1\}$, $j\in\{1,\ldots,r_2\}$, and $a_i,b_j\in\{0,1\}$. Let $F^1$ (resp. $F^0$) be the direct summand of $M$ generated by all those $e_i$ and $f_j$ with $a_i=1$ and $b_j=1$ (resp. with $a_i=0$ and $b_j=0$). Let $\{e_i\otimes f_j^*|i\in\{1,\ldots,r_1\}\;,\;j\in\{1,\ldots,r_2\}\}$ be the natural $W(k)$-basis for $\Hom(M_2,M_1)=M_1\otimes_{W(k)} \Hom(M_2,W(k))$. It is easy to see that $E$ is $W(k)$-generated by $p^{c_{ij}}e_i\otimes f_j^*$ for suitable $c_{i,j}\in\dbN\cup\{0\}$. We have $a_i-b_j=0$ (resp. $a_i-b_j=1$) if and only if  $p^{c_{ij}}e_i\otimes f_j^*\in F^0(E)$ (resp. $p^{c_{ij}}e_i\otimes f_j^*\in F^1(E)$). Let $\Lambda\subseteq \{1,\ldots,r_1\}\times \{1,\ldots,r_2\}$ be such that the following property holds: If $(i,j)\in\Lambda$, then we have $a_i-b_j\in\{0,1\}$ and $(\pi_1(i),\pi_2(j))\in\Lambda$. This implies that for each $(i,j)\in\Lambda$ we have $c_{i,j}=0$ and thus $\{e_i\otimes f_j^*|(i,j)\in\Lambda\}$ is a $W(k)$-basis for a direct summand $E_0$ of either $E$ or $V_-$ for which the axioms 1.4 (i) and (ii) hold. Moreover, if $D_{W(k)}$ is the $p$-divisible group over $W(k)$ which lifts $D$ and whose  filtered Dieudonn\'e module is $(M,F^1,\phi)$, then the axioms 1.4 (iii) and (iv) hold for the pair $(E_0,D_{W(k)})$. From Basic Theorem B we get that to the pair $(E_0,D_{W(k)})$ corresponds a smooth formal subscheme $\scrG^{E_0}$ of $\scrD$ that has a canonical structure of a formal Lie group over $\Spf(W(k))$ isomorphic to the one associated to $\dbD^{-1}(E)_{W(k)}^{\text{t}}$ of property 1.6 (i) and that has an origin which corresponds to $D_{W(k)}$. 
Note that $\scrG_k^{E_0}=\scrG^{E_0}\cap \scrD_k$ is a formal Lie subgroup of $\scrG_k^E$, cf. Proposition 3.6 (a). 

Numerical example: If $r_1=r_2=3$, $\alpha_1={1\over 3}$, $\alpha_2={2\over 3}$, $\pi_1=\pi_2=(1\,2\,3)$, $a_1=b_1=b_2=1$, and $a_2=a_3=b_3=0$, then we can take $\Lambda=\{(1,1),(2,2),(3,3),(1,2),(2,3),(3,1)\}$; the ranks of $E$, $E_0$, and $F^0(E_0)$ are $9$, $6$, and $4$ (respectively).

\medskip\smallskip\noindent
{\bf 1.8. Example.} Let $M=F^1\oplus F^0$, $\mu:\dbG_m\to \pmb{\GL}_M$, and $D_{W(k)}$ be as in Subsubsection 1.3.1 and Subsection 1.4. Suppose that there exist a reductive, closed subgroup scheme $G$ of $\pmb{\GL}_M$, a maximal torus $T$ of $G$, and two opposite Borel subgroup schemes $B_+$ and $B_-$ of $G$ that contain $T$, such that the following five properties hold:

\medskip
{\bf (i)} the four Lie algebras $\Lie(G)[{1\over p}]$, $\Lie(B_+)[{1\over p}]$, $\Lie(B_-)[{1\over p}]$, and $\Lie(T)[{1\over p}]$ are stable under $\phi$;

\smallskip
{\bf (ii)} the cocharacter $\mu:\dbG_m\to \pmb{\GL}_M$ factors through $T$;

\smallskip
{\bf (iii)} all Newton polygon slopes of $(\Lie(B_+)[{1\over p}],\phi)$ (resp. of $(\Lie(B_-)[{1\over p}],\phi)$) are non-negative (resp. are non-positive);

\smallskip
{\bf (iv)} the Borel subgroup scheme $B_+$ (resp. $B_-$) normalizes $F^1$ (resp. $F^0$);

\smallskip
{\bf (v)} if $E:=V_-\cap\Lie(G)$, then we have $E^2=0$.

\medskip
The $\sigma$-linear automorphism $\sigma_{\phi}$ of $M$ (and thus also of $\End(M)$) normalizes $\Lie(G)$, $\Lie(B^+)$, $\Lie(B_-)$, and $\Lie(T)$ (cf. (i) and (ii)). Based on this and (iii) and (iv), we have $F^{-1}(E)=\Lie(B^-)\cap \Hom(F^1,F^0)$ and $E=F^0(E)\oplus F^{-1}(E)$ is a direct summand of $\Lie(B_-)$. More precisely, $F^0(E)$ is the maximal direct summand of $\Lie(B_-)$ that has the following three properties: (a) it is normalized by all conjugates of $\mu_L$ under integral powers of $\sigma_{\phi}$, (b) it has no element fixed by all these conjugates of $\mu_L$, and (c) it is fixed by $\mu_L$. As $E=F^0(E)\oplus F^{-1}(E)$, the axiom 1.4 (iii) holds for $(E,D_{W(k)})$ and $\sigma_{\phi}=\phi\mu(p)$ normalizes $E[{1\over p}]\cap\Lie(B_-)=E$. We have $\phi(F^0(E))\oplus p\phi(F^{-1}(E))=\sigma_{\phi}(E)=E$ i.e., the axiom 1.4 (iv) holds for $(E,D_{W(k)})$. From this and the fact that $E$ is a direct summand of $\End(M)$ and therefore also of $V_-$, we get that the axiom 1.4 (i) holds. Axiom 1.4 (ii) is implied by (v). Thus the axioms 1.4 (i) to (iv) hold for $E$ and $(E,D_{W(k)})$. 

If the derived group scheme $G^{\der}$ of $G$ is $\pmb{\text{Sp}}_{2n}^s$ (resp. $\pmb{\text{Spin}}_{2n+1}^s$) for some $s\in\dbN$ and if the representation of $G^{\der}$ on $M$ is a direct sum of irreducible representations of rank $2n$ (resp. $2^n$), then the property (v) is implied by the properties (i) to (iv). The triple $(M,\phi,G)$ is called an {\it ordinary $F$-crystal ($p$-divisible object) with a reductive group} over $k$ and $F^1$ is called the {\it canonical lift} of $(M,\phi,G)$, cf. [Va2, Def. 4.3.3]. Based on the Basic Theorem B, to $(E,D_{W(k)})$ corresponds a formal subscheme $\scrG^E$ of $\scrD$ that has a natural structure of a formal Lie group over $\Spf(W(k))$ isomorphic to the formal Lie group over $\Spf(W(k))$ associated to $\dbD^{-1}(E)_{W(k)}^{\text{t}}$. As the lift $F^1$ of $(M,\phi,G)$ is canonical (cf. [Va2, Prop. 4.3.2 (b)]), for $p>2$ (resp. for $p=2$) there exists a canonical choice (resp. a finite number of quasi-canonical choices) for $\dbD^{-1}(E)_{W(k)}$ and this (resp. and often a suitable such) choice makes the formal Lie group structure on $\scrG^E$ to be canonical.

\medskip\smallskip\noindent
{\bf 1.9. Example.} Suppose $\scrW$ has four elements $\alpha_1<\alpha_2<\alpha_3<\alpha_4$. Let $E$ be the maximal $W(k)$-submodule of $V_-$ such that $(E,p\phi)$ is a Dieudonn\'e module and $E$ is contained in the direct summand 
$\Hom(W(\alpha_4),W(\alpha_3))\oplus \Hom(W(\alpha_4),W(\alpha_1))\oplus \Hom(W(\alpha_2),W(\alpha_1))$
of $\End(M[{1\over p}])=\oplus_{i,j=1}^4 \Hom(W(\alpha_i),W(\alpha_j))$. 
As $E^2=0$, the axioms 1.4 (i) and (ii) hold for $E$. There exists a short exact sequence $0\to D_2\to D\to D_1\to 0$ such that $\scrG_k^E$ is a formal Lie subgroup of $Ext^1(D_1,D_2)$, cf. Corollary 3.7. The choice of $E$ forces the Newton polygon slopes of $D_1$ (resp. $D_2$) to be $\alpha_1$ and $\alpha_3$ (resp. $\alpha_2$ and $\alpha_4)$. 

This example can be adapted to the limit case when $\alpha_2=\alpha_3$ (i.e., when $\scrW=\{\alpha_1,\alpha_2,\alpha_4\}$ and we have a direct sum decomposition $(W(\alpha_2),\phi)=(\tilde W(\alpha_2),\phi)\oplus (\tilde W(\alpha_3),\phi)$ into non-zero $F$-crystals). In this case, $\alpha_2$ is a Newton polygon slope of both $D_1$ and $D_2$.

\medskip\smallskip\noindent
{\bf 1.10. Refined structures.} Let $\scrC_k$ be the largest smooth formal subscheme of $\scrD_k$ such that the restriction of $\grD^{\text{f}}$ to it is geometrically constant (see Section 5). Referring to the Basic Theorem A, there exist many operations between different $\scrG_k^E$'s or between $\scrC_k$ and some $\scrG_k^E$'s. To be short, here we will mention briefly only two such operations (to be compared with [Va4, Chs. 5 and 9], where these operations show up first). If $E_1$, $E_2$, and $E_3$ are $W(k)$-submodules of $V_-$ for which the axioms 1.4 (i) and (ii) hold and moreover we have $[E_1,E_2]\subseteq E_3$, then one can check that there exists a natural {\it commutator bilinear formal morphism} $\scrG_k^{E_1}\times_{\Spf(k)} \scrG_k^{E_2}\to \scrG_k^{E_3}$. If $[E_1,V_-]=0$, then one can check that we have a {\it formal group action} $\scrG_k^E\times_{\Spf(k)} \scrC_k\to \scrC_k$. 

Based on these operations, one can canonically associate to $D$ a {\it nilpotent formal Lie group} $\scrN_k$ over $\Spf(k)$. We state here what we call the Main Problem of the formal Lie group theory for formal deformation spaces of $p$-divisible groups over $k$. 

\medskip\noindent
{\bf Main Problem.} Let $\scrC_k$ be as above. 

\medskip
{\bf (a)} Show that for each $p$-divisible group $D$ (resp. each $D$ that is a direct sum of isoclinic $p$-divisible groups) over $k$, $\scrC_k$ has a canonical structure of a nilpotent formal Lie group over $\Spf(k)$ which is isomorphic to $\scrN_k$ and for which the following properties hold:

\medskip\noindent
{\bf (a.i)} each $\scrG_k^E$ as in the Basic Theorem A is a formal Lie subgroup of it;

\smallskip\noindent
{\bf (a.ii)} it produces naturally the above operations (i.e., the mentioned commutator bilinear formal morphisms, formal group actions, etc.);

\smallskip\noindent
{\bf (a.iii)} it has a composite series whose factors are formal Lie groups associated to $p$-divisible groups over $k$ isogenous to (resp. isomorphic to) the $Z_{+,-}(\{(\alpha,\beta)\})$'s of Definition 4.1.4. 

\medskip
{\bf (b)} Classify all local Traverso strata $\scrC_k(G,\scrB_G)_{\text{red}}$ of Definition 5.2.1 that are formal Lie subgroups of $\scrC_k$ (see Example 5.2.4).

\smallskip
{\bf (c)} Classify all formal subschemes of $\scrD_k$ (resp. of $\scrD$) that have a canonical (resp. natural) structure of a nilpotent formal Lie group over $\Spf(k)$ (resp. over $\Spf(W(k))$).  
 
\medskip
The Basic Theorem A and Proposition 3.6 (resp. [Ch2]) essentially solve the Main Problem for the case when the set $\scrW$ has two elements (resp. for the case when $D$ is a direct sum of two isoclinic $p$-divisible groups of distinct Newton polygon slopes). Certainly, even the non-canonical forms of (a) are much stronger forms than the previously suggested approaches to the study of $\scrC_k$ based on the fact that $\scrC_k\times_{\scrD_k} \grD^{\text{f}}$ has an increasing Newton polygon filtration. [It allows a filtering of $\scrC_k$ in such a way that the fibres of the resulting quotient formal morphisms are torsors under suitable formal Lie groups, etc.]. 

\medskip\smallskip\noindent
{\bf 1.11. On contents.} In Section 2 we present complements on Faltings deformation theory (see [Fa, Sect. 7]) that are needed in the paper. The proofs of the Basic Theorems A and B are carried on in Section 3. They rely on Section 2 and on the methods used in [Va4, Ch. 9]; they are short, elementary, and foundational. In Section 4 we extend the concepts of {\it sign $p$-divisible groups} and of {\it duality properties} first mentioned in [Va4, Ch. 5]. Subsection 5.1 contains a direct description of the smooth formal subscheme $\scrC_k$ of $\scrD_k$. Subsections 5.2 and 5.3 define abstractly the most general type of local Traverso strata. Example 5.4 pertains to the principally quasi-polarized context. The above notations $p$, $k$, $W(k)$, $B(k)$, $\sigma=\sigma_k$, $c$, $d$, $r=c+d$, $D$, $(M,\phi)$, $R=W(k)[[x_1,\ldots,x_{cd}]]$, $\scrD=\Spf(R)$, $\scrD_k=\Spf(R/pR)$, $\grD$, $\grD^{\text{f}}$, $M[{1\over p}]=\oplus_{\alpha\in\scrW} W(\alpha)$, $\End(M[{1\over p}])=L_+\oplus L_0\oplus L_-$, $V_+$, $V_-$, $\bar M$, $\bar F^1$, $\bar\phi$, $\scrT$, $\nu:\End(M)\twoheadrightarrow \scrT$, $M=F^1\oplus F^0$, $\mu:\dbG_m\to \pmb{\GL}_M$, $\mu_L:\dbG_m\to \pmb{\GL}_{\End(M)}$, and $\sigma_{\phi}$ are used throughout the paper. If the pair $(E,p\phi)$ is a Dieudonn\'e module over $k$, let $\dbD^{-1}(E)$ be a $p$-divisible group over $k$ whose Dieudonn\'e module is isomorphic to $(E,p\phi)$.

\bigskip
\noindent
{\boldsectionfont 2. Complements on Faltings deformation theory}
\bigskip 

Let $n\in\dbN\cup\{0\}$. Let $S:=W(k)[[x_1,\ldots,x_n]]$. Let $I$ be the maximal ideal of $S$. Let $J:=(x_1,\ldots,x_n)$; it is an ideal of $S$. Let $\Phi_S$ be the Frobenius lift of $S$ that is compatible with $\sigma$ and that takes $x_i$ into $x_i^p$ for all $i\in\{1,\ldots,n\}$. We identify $\Omega_S^\wedge:=\oplus_{i=1}^n Sdx_i$ with the $p$-adic completion of the $S$-module of differentials $\Omega_S$. Let $d\Phi_S:\Omega_S^\wedge\to\Omega_S^\wedge$ be the differential map of $\Phi_S$. An {\it $F$-crystal} over $S/pS$ is a triple of the form $(Q,\Phi_{Q},\Delta)$, where $Q$ is a free $S$-module of finite rank, where $\Phi_{Q}:Q\hookrightarrow Q$ is a $\Phi_S$-linear monomorphism, and where $\Delta:Q\to Q\otimes_S \Omega_S^\wedge$ is an integrable and topologically nilpotent {\it connection} on $Q$, such that the following two properties hold:

\medskip
-- the $S$-module which is the quotient of $Q$ by the $S$-span of $\im(\Phi_{Q})$, is annihilated by some non-negative integral power of $p$;

\smallskip
-- we have an identity $\Delta\circ\Phi_{Q}=(\Phi_{Q}\otimes d\Phi_S)\circ\Delta$. 

\medskip\noindent
If $F^1_{Q}$ is a direct summand of $Q$ such that $\Phi_{Q}$ induces an $S$-linear isomorphism $({1\over p}F^1_{Q}+Q)\otimes_S {}_{\Phi_S} S\arrowsim Q$, then we call the quadruple $(Q,F^1_{Q},\Phi_Q,\Delta)$ a {\it filtered $F$-crystal} over $S/pS$.

Faltings deformation theory for such filtered $F$-crystals $(Q,F^1_{Q},\Phi_Q,\Delta)$  over $S/pS$, can be summarized by the following two statements:

\medskip
$\bullet$ The connection $\Delta$ is uniquely determined by the identity $\Delta\circ\Phi_{Q}=(\Phi_{Q}\otimes d\Phi_S)\circ\Delta$ and it is automatically integral and topologically nilpotent.

\smallskip
$\bullet$ If $\Gamma$ is a $p$-divisible group over $\Spf(W(k))$ whose filtered Dieudonn\'e module is $(Q,F^1_Q,\Phi_Q)$ modulo $J$ (it exists even if $p=2$ as $k$ is algebraically closed), then the filtered $F$-crystal $(Q,F^1_Q,\Phi_Q,\Delta)$ over $S/pS$ is associated (via the crystalline Dieudonn\'e functor $\dbD$) to a $p$-divisible group $\grG$ over $\Spf(S)$ that lifts $\Gamma$.

\medskip
For a general $\Phi_Q$, it is hard to either (i) compute explicitly the connection $\Delta$ or (ii) decide if the restriction of $\grG$ to some smooth formal subscheme of $\Spf(S/pS)$ is or is not geometrically constant.

Below we will introduce a setting (see Subsection 2.1) that involves an explicit $Q$ and a particular type of $\Phi_Q$ for which we can compute the connection $\Delta$ (see Proposition 2.6) and for which we can later on find a large smooth formal subscheme of $\Spf(S/pS)$ over which $\grG$ restricts to a geometrically constant $p$-divisible group (see Lemma 3.1). The computation of $\Delta$ plays a role in the identification inside $\End(Q)$ of free $S$-modules that are left invariant by $\Delta$ and that give birth naturally to (filtered) $F$-crystals over $S/pS$. 

\medskip\smallskip\noindent
{\bf 2.1. A setting.} Let $\scrB:=(v_1,\ldots,v_n)$ be an $n$-tuple of elements of $\End(M)$ for which the following property holds:

{\it If $N:=\sum_{i=1}^n W(k)v_i\subseteq\End(M)$ and $n_0:=\dim_k(\nu(N))\in\dbN\cup\{0\}$, then the set $\{\nu(v_i)|i\in\{1,\ldots,n_0\}\}$ is a $k$-basis for $\nu(N)$.}

Let $S_0:=W(k)[[x_1,\ldots,x_{n_0}]]$. We identify  $\Spec(S_0)$ with the closed subscheme of $\Spec(S)$ defined by the ideal $(x_{n_0+1},\ldots,x_n)$ of $S$. Let $N_0:=\sum_{i=1}^{n_0} W(k)v_i$. As the map $N_0/pN_0\to\scrT$ induced by $\nu$ is injective, $N_0$ has rank $n_0$ and it is a direct summand of $\End(M)$ and therefore also of $N$.  

Let $U$ be the smooth, affine scheme over $\Spec(W(k))$ defined by the following rule on valued points: if $A$ is a commutative $W(k)$-algebra, then 
$$U(A):=\{x\in 1_{M\otimes_{W(k)} A}+N\otimes_{W(k)} A|\text{det}(x)\,\,\text{is}\,\,\text{an}\,\,\text{invertible}\,\,\text{element}\,\,\text{of}\,\,A\}$$ 
is the subset of $1_{M\otimes_{W(k)} A}+N\otimes_{W(k)} A$ formed by those elements whose images in $\break\End(M)\otimes_{W(k)} A$ belong to $\pmb{\GL}_M(A)$. We have a natural morphism $h_U:U\to \pmb{\GL}_M$ whose fibre over $\Spec(B(k))$ is a closed embedding. Thus if $A$ is an integral domain of characteristic $0$, then we have an injective map $h_U(A):U(A)\hookrightarrow \pmb{\GL}_M(A)$; to ease notations, we always write $U(A)\subseteq \pmb{\GL}_M(A)$ and we do not use $h_U$ besides this and the next paragraph. Let $U_0$ be the smooth, closed subscheme of $U$ defined by the following rule $U_0(A):=(1_{M\otimes_{W(k)} A}+N_0\otimes_{W(k)} A)\cap U(A)$.

If $N$ is a direct summand of $\End(M)$, then $h_U$ is a closed embedding. If $N$ is an algebra (i.e., we have $N^2\subseteq N$), then the inverse of $x\in U(A)$ belongs to $A[x]$ (cf. Cayley--Hamilton theorem) and thus also to $U(A)$; this implies that:

\medskip
{\bf (i)} $U$ is a group scheme over $\Spec(W(k))$ whose Lie algebra is $\Lie(N)$, and 

\smallskip
{\bf (ii)} $h_U$ is a homomorphism of group schemes over $\Spec(W(k))$. 

\medskip
If $U$ is a group scheme, then any $W(k)$-submodule of $\End(M)$ which is normalized by $U$ under inner conjugation, will be called a left $U$-module. The element
$$u_{\text{univ}}:=1_{M\otimes_{W(k)} S}+\sum_{i=1}^n v_i\otimes x_i\in U(S)\subseteq \pmb{\GL}_M(S)$$ is universal in the sense that the natural morphism $u_{\text{univ}}:\Spec(S)\to U$ is formally smooth (formally \'etale, if $N$ has rank $n$). We consider the $\Phi_S$-linear monomorphism 
$$\Phi_N:=u_{\text{univ}}(\phi\otimes\Phi_S):M\otimes_{W(k)} S\hookrightarrow M\otimes_{W(k)} S.$$  
\indent
If $\ast$ is a free $W(k)$-module of finite rank, let $\delta_0$ be the flat connection on $\ast\otimes_{W(k)} S$ that annihilates $\ast\otimes 1$. Let $D_{W(k)}$ be an arbitrary $p$-divisible group over $W(k)$ that lifts $D$. Let $F^1$ be the direct summand of $M$ which is the Hodge filtration of $D_{W(k)}$. 
The reduction of $\Phi_N$ modulo $J$ is $\phi$ and $\Phi_N$ induces naturally an $S$-linear isomorphism $\break ({1\over p}F^1+M)\otimes_{W(k)} S\otimes_S {}_{\Phi_S} S\arrowsim M\otimes_{W(k)} S$.
From this and [Fa, Sect. 7, Thm. 10] we get:

\medskip\smallskip\noindent
{\bf 2.2. Theorem.} {\it {\bf (a)} There exists a unique connection 
$$\nabla:M\otimes_{W(k)} S\to M\otimes_{W(k)} S\otimes_S \Omega_S^\wedge=M\otimes_{W(k)} \Omega_S^\wedge$$ 
such that the $\Phi_S$-linear endomorphism $\Phi_N:M\otimes_{W(k)} S\to M\otimes_{W(k)} S$ is horizontal with respect to $\nabla$ i.e., we have an identity 
$$\nabla\circ \Phi_N=(\Phi_N\otimes d\Phi_S)\circ\nabla:({1\over p}F^1+M)\otimes_{W(k)} S\to M\otimes_{W(k)} \Omega_S^\wedge.$$ 
Moreover the connection $\nabla$ is integrable and topologically nilpotent. 

\smallskip
{\bf (b)} There exists a $p$-divisible group $\grH^{\text{f}}$ over $\Spf(S)$ which modulo $J$ is the $p$-divisible group over $\Spf(W(k))$ defined by $D_{W(k)}$ and whose filtered $F$-crystal over $S/pS$ is the quadruple $(M\otimes_{W(k)} S,F^1\otimes_{W(k)} S,\Phi_N,\nabla)$.}

\medskip
Let $\grH$ be the $p$-divisible group over $\Spec(S)$ defined by $\grH^{\text{f}}$. As $\{x_1,\ldots,x_n\}$ is a $p$-basis (in the sense of [BM, Def. 1.1.1]) for the regular ring $S/pS$, from [BM, Thm. 4.1.1] we get that $\grH_{S/pS}$ is (up to isomorphism) the unique $p$-divisible group over $\Spec(S/pS)$ whose $F$-crystal over $S/pS$ is $(M\otimes_{W(k)} S,\Phi_N,\nabla)$. For all $m\in\dbN$, the ideal $pJ$ of $S$ has a nilpotent divided power structure modulo $pJ^m$. From the last two sentences and the Grothendieck--Messing deformation theory (see [Me, Chs. 4 and 5]), we get that $\grH^{\text{f}}$ and $\grH$ are uniquely determined (up to unique isomorphisms) by (b). Let 
$$q:\Spec(S)\to\Spec(R)$$ 
be the unique morphism such that we have an identification $\grH=q^*(\grD)$ that lifts the identity automorphism of $D$. Let $q_{k}:\Spec(S/pS)\to\Spec(R/pR)$, $q_0:\Spec(S_0)\to\Spec(R)$, and $q_{0k}:\Spec(S_0/pS_0)\to\Spec(R/pR)$ be the natural restriction morphisms defined by $q$. We consider the $S$-linear map
$$\grK:\oplus_{i=1}^n S{{\partial}\over {\partial x_i}}\to \Hom(F^1,M/F^1)\otimes_{W(k)} S$$
that is the Kodaira--Spencer map of $\nabla$. Let
$$\bar\grK:\oplus_{i=1}^n k{\partial\over {\partial x_i}}\to \Hom(F^1,M/F^1)\otimes_{W(k)} k=\Hom(\bar F^1,\bar M/\bar F^1)=\scrT$$ 
be the reduction of $\grK$ modulo $I$. Example: if $\scrB=(0,\ldots,0)$, then $\nabla=\delta_0$ and thus $\grK=0$.

\medskip\noindent
{\bf 2.2.1. Definition.} By the {\it deformation package} of $\scrB=(v_1,\ldots,v_n)\in\End(M)^{\times n}$, we mean the following $14$-tuple (whose entries were introduced above):
$$\dbP(\scrB):=(n,N,n_0,N_0,U,U_0,u_{\text{univ}},\Phi_N,\nabla,\grH_{S/pS},q_k,q_{0k},\grK,\bar\grK).$$ 
\noindent
{\bf 2.3. Lemma.} {\it We have $\im(\bar\grK)=\nu(N)=\nu(N_0)$; thus $\dim(\im(\bar\grK))=n_0$. The morphism $q_0:\Spec(S_0)\hookrightarrow\Spec(R)$ is a closed embedding and the tangent space of the closed subscheme $\Spec(S_0/pS_0)$ of $\Spec(R/pR)$ is equal to $\nu(N_0)$.}

\medskip
\proof
Let $\delta_1$ be the flat connection on $M\otimes_{W(k)} S$ that annihilates the free $W(k)$-module $u_{\text{univ}}(M\otimes 1)$; we have $\delta_1=\delta_0-u_{\text{univ}}^{-1}du_{\text{univ}}$. As $d(\Phi_S(x_i))=d(x_i^p)=px_i^{p-1}\in I$, the reduction modulo $I$ of $(\Phi_N\otimes d\Phi_S)\circ\nabla$ is trivial. Due to this and the equation $\nabla\circ \Phi_N=(\Phi_N\otimes d\Phi_S)\circ\nabla$, the reductions modulo $I$ of $\delta_1$ and $\nabla$ coincide. As $u_{\text{univ}}^{-1}$ modulo $I$ is $1_{\bar M}$ and as $du_{\text{univ}}=\sum_{i=1}^n n_i\otimes dx_i$, we get that $\im(\bar\grK)=\sum_{i=1}^n k\nu(n_i)=\nu(N)=\nu(N_0)$. Thus $\dim_k(\im(\bar\grK))=n_0$. The restriction of $\bar\grK$ to $\oplus_{i=1}^{n_0} k{\partial\over {\partial x_i}}$ is injective and its image is $\nu(N_0)$. Thus the tangent map of $q_0$ (at $k$-valued points) is an injective map whose image is $\nu(N_0)$. This implies that $q_{0k}$ is a closed embedding. Thus $q_0$ itself is a closed embedding.\endproof

\medskip\smallskip\noindent
{\bf 2.4. Fact.} {\it Let $K$ be an algebraically closed field that contains $k$. Let $h:\Spec(K)\to\Spec(S)$ be a morphism. Then the Dieudonn\'e module of $h^*(\grH)$ is isomorphic to $(M\otimes_{W(k)} W(K),u_h(\phi\otimes\sigma_K))$ for some element $u_h\in U(W(K))$.} 

\medskip
\proof
Let $T_h:\Spec(W(K))\to\Spec(S)$ be the Teichm\"uller lift of $h$ (see [Ka1, Subsect. 2.4, p. 145]). Let $u_h\in U(W(K))\leqslant \pmb{\GL}_M(W(K))$ be the composite of $T_h\in\Spec(S)(W(K))$ with the point $u_{\text{univ}}\in U(S)$. As $T_h$ is a Teichm\"uller lift of $h$, the Dieudonn\'e module of $h^*(\grH)$ is naturally identified with $(M\otimes_{W(k)} W(K),u_h(\phi\otimes\sigma_K))$.\endproof

\medskip\smallskip\noindent
{\bf 2.5. Fact.} {\it Suppose that there exists a $W(k)$-submodule $E$ of $M$ with $E^2=0$ and $N\subseteq E$ and that there exists an element $t\in\End(M[{1\over p}])$ such that $[E,t]=E$. Let $\tilde E:=E\oplus W(k)t$; it is a Lie algebra over $W(k)$ contained in $\End(M[{1\over p}])$. Then $\tilde E$ is a left $U$-module. If $N$ is a direct summand of $E$, then the homomorphism $U\to \pmb{\GL}_{\tilde E}$ is a closed embedding.}

\medskip
\proof
Let $x\in N$. As $E^2=0$ and $[x,t]\in E$, we have $xtx=xxt-x[x,t]=0$. Thus
$$(1_M+x)t(1_M+x)^{-1}=(1_M+x)t(1_M-x)=t+[x,t]-xtx=t+[x,t]\in\tilde E\otimes_{W(k)} A.\leqno (1)$$ 
As $N\subseteq E$ and $E^2=0$, we have $N^2=0$ and therefore $U$ is a group scheme over $\Spec(W(k))$ that fixes $E$ under inner conjugation. From the last two sentences we get that $\tilde E$ is a left $U$-module.  Let $\tilde U$ be the smooth, commutative, closed subgroup scheme of $\pmb{\GL}_{\tilde E}$ defined by the following rule on valued points: if $A$ is a commutative $W(k)$-algebra, then we have $\tilde U(A):=1_{\tilde E\otimes_{W(k)} A}+\Hom(\tilde E/E,E)\otimes_{W(k)} A\leqslant \pmb{\GL}_{\tilde E}(A)$. As both $E$ and $\tilde E/E$ are trivial left $U$-modules, the homomorphism $U\to \pmb{\GL}_{\tilde E}$ factors through $\tilde U$. The differential map of the resulting homomorphism $l_U:U\to\tilde U$ is a restriction of the adjoint representation of $\Lie(\pmb{\GL}_{M[{1\over p}]})$; thus $dl_U:\Lie(U)\to\Lie(\tilde U)$ is defined by the rule:
$$\text{if}\,\,x\in\Lie(U)=N,\,\,\text{then}\,\,\forall y\in\tilde E\,\,\text{we}\,\,\text{have}\,\,dl_U(x)(y)=[x,y]\in E.\leqno (2)$$

We assume that $N$ is a direct summand of $E$. Based on Formula (2) and the identity $[E,t]=E$, it is easy to see that $dl_U$ is an injective $W(k)$-linear map whose image is a direct summand of $\Lie(\tilde U)=\Hom(\tilde E/E,E)$. 

But $l_U$ is uniquely determined by its generic fibre and thus (cf. [Bo, Ch. II, Subsect. 7.1]) by $dl_U[{1\over p}]$. As $U$ and $\tilde U$ are vector group scheme over $\Spec(W(k))$ (cf. their definitions), we get that for each $A$ as above, the homomorphism $l_U(A):U(A)\to\tilde U(A)$ maps $1_{M\otimes_{W(k)} A}+y$ with $y\in N\otimes_{W(k)} A=\Lie(U)\otimes_{W(k)} A$ to $1_{\tilde E\otimes_{W(k)} A}+dl_U(y)\in 1_{\tilde E\otimes_{W(k)} A}+\Lie(\tilde U)\otimes_{W(k)} A=\tilde U(A)$ (to be compared with Formulas (1) and (2)). From this and the previous paragraph we get that $l_U$ is a closed embedding.\endproof

\medskip\smallskip\noindent
{\bf 2.6. Proposition.} {\it Let the connection $\nabla$ on $M\otimes_{W(k)} S$ be as in the Theorem 1.2 (a). Let $E$ be a $W(k)$-submodule of $V_-$ such that $p\phi(E)\subseteq E$ and $E^2=0$. Suppose that $N\subseteq E$.

\medskip
{\bf (a)} Then there exists $\omega\in E\otimes_{W(k)} \Omega_S^\wedge$ such that we have $\nabla=\delta_0+\omega$.

\smallskip
{\bf (b)} Let $t\in\End(M[{1\over p}])$ be such that $[E,t]=E$. Let $\tilde E:=E\oplus W(k)t$. The connection $\nabla$ on $M\otimes_{W(k)} S$ induces a connection on $\End(M[{1\over p}])\otimes_{B(k)} S[{1\over p}]$ that restricts to a connection $\tilde\nabla=\delta_0+\tilde\omega$ on $\tilde E\otimes_{W(k)} S$, where $\tilde\omega\in\Hom(\tilde E/E,E)\otimes_{W(k)} \Omega_S^\wedge\subseteq \End(\tilde E)\otimes_{W(k)} \Omega_S^\wedge$.}

\medskip
\proof
Let $\{c_1,\ldots,c_r\}$ be a $W(k)$-basis for ${1\over p}F^1+M$. Let $m$ be the rank of $E$ and let $\{e_1,\ldots,e_m\}$ be a $W(k)$-basis for $E$. We write $p\phi(e_l)=\sum_{j=1}^m a_{jl}e_j$ and $v_i=\sum_{j=1}^m b_{ji}e_j$, where $(a_{jl})_{1\le j,l\le m}\in M_{m\times m}(W(k))$ and $(b_{ji})_{1\le j\le m,1\le i\le n}\in M_{m\times n}(W(k))$. 

Each element of $E\otimes_{W(k)} \Omega_S^\wedge$ can be written as
$$\omega=\sum_{i=1}^n \sum_{l=1}^m e_l\otimes w_{l,i}dx_i,\;\;\;\text{where}\;\;\; w_{l,i}\in S.$$ 
Based on the uniqueness part of Theorem 2.2 (a), to prove (a) it suffices to show that there exists functions $w_{l,i}\in S$ such that for all $s\in\{1,\ldots,r\}$ we have an identity
$$(\delta_0+\omega)(\Phi_N(c_s\otimes 1))=(\Phi_N\otimes d\Phi_S)((\delta_0+\omega)(c_s\otimes 1)).\leqno (3)$$
As $E^2=0$, we have $e_lv_i=0$ and $v_i\phi(e_l)=0$. Based on this, the left hand side of (3) is 
$$(\delta_0+\omega)[\phi(c_s)\otimes 1+\sum_{i=1}^n v_i(\phi(c_s))\otimes x_i]=\sum_{i=1}^n v_i(\phi(c_s))\otimes dx_i+\sum_{i=1}^n\sum_{j=1}^m e_j(\phi(c_s))\otimes w_{j,i}dx_i$$
$$=\sum_{j=1}^m e_j(\phi(c_s))\otimes [\sum_{i=1}^n(b_{ji}+w_{j,i})dx_i]$$ 
and the right hand side of (3) is 
$$(\Phi_N\otimes d\Phi_S)(\sum_{i=1}^n \sum_{l=1}^m e_l(c_s)\otimes w_{l,i}dx_i)=\sum_{i=1}^n \sum_{l=1}^m u_{\text{univ}}(\phi(e_l(c_s)))\otimes p\Phi_S(w_{l,i})x_i^{p-1}dx_i$$
$$=\sum_{i=1}^n \sum_{l=1}^m p\phi(e_l)(\phi(c_s))\otimes \Phi_S(w_{l,i})x_i^{p-1}dx_i=\sum_{i=1}^n \sum_{l=1}^m\sum_{j=1}^m a_{jl}e_j(\phi(c_s))\otimes \Phi_S(w_{l,i})x_i^{p-1}dx_i$$
$$=\sum_{j=1}^m e_j(\phi(c_s))\otimes [\sum_{i=1}^n (\sum_{l=1}^m a_{jl}\Phi_S(w_{l,i}))x_i^{p-1}dx_i].$$
Thus (3) holds if and only if for all $i\in\{1,\ldots,n\}$, we have
$$b_{ji}+w_{j,i}=\sum_{l=1}^m a_{jl}\Phi_S(w_{l,i})x_i^{p-1}\;\;\;\;\forall j\in\{1,\ldots,m\}.\leqno (4)$$
Defining $b_i:=(b_{1i},\ldots,b_{mi})$, $w_i:=(w_{1,i},\ldots,w_{m,i})$, and $\grO_i(w_i):=(\grO_{1,i}(w_i),\ldots,\grO_{m,i}(w_i))$ with $\grO_{j,i}(w_i):=\sum_{l=1}^m a_{jl}\Phi_S(w_{l,i})x_i^{p-1}$, for $i\in\{1,\ldots,n\}$ the system (4) becomes
$$b_i+w_i=\grO_i(w_i).\leqno (5)$$
As the operator $\grO_i:S^m\to S^m$ is topologically nilpotent in the $I$-adic topology, the system (5) has the unique solution 
$$w_i=\sum_{w=0}^{\infty} \grO_i^w(-b_i).$$ 
Thus there exists a unique $mn$-tuple $(w_{1,1},\ldots,w_{m,n})\in S^{mn}$ such that the equations (3) hold for all $s\in\{1,\ldots,r\}$. This ends the proof of (a).

We check (b). For $x\in\End(M)$, the connection $\nabla$ on $\End(M[{1\over p}])\otimes_{B(k)} S[{1\over p}]$ maps $x\otimes 1$ to $\omega(x\otimes 1)=\sum_{i=1}^n \sum_{l=1}^m [e_l,x]\otimes w_{l,i}dx_i$. As $[E,E]=0$ and $[E,\tilde E]\subseteq E$, $\nabla$ induces a connection $\tilde\nabla$ on $\tilde E\otimes_{W(k)} S$ that restricts to the flat connection $\delta_0$ on $E\otimes_{W(k)} S$ and that induces the flat connection $\delta_0$ on $\tilde E/E\otimes_{W(k)} S$. Therefore we can write $\tilde\nabla=\delta_0+\tilde\omega$, where $\tilde\omega\in\Hom(\tilde E/E,E)\otimes_{W(k)} S\subseteq \End(\tilde E)\otimes_{W(k)} S$. Thus (b) holds.\endproof

\medskip\smallskip\noindent
{\bf 2.7. Corollary.} {\it Let $E$ be a $W(k)$-submodule of $V_-$ such that the axioms 1.4 (i) and (ii) hold. Let $\scrB_E$ be the set of $n$-tuples $\tilde\scrB$ of elements of $E$ which are ordered $W(k)$-bases for direct summands $\tilde N$ of $E$ that have the following two properties: (i) $\nu(\tilde N)=\nu(E)$, and (ii) the $k$-linear map $\tilde N/p\tilde N\to \scrT$ induced by $\nu$ is injective. Suppose that $\scrB\in\scrB_E$; thus $N\subseteq E$. Then the closed embedding $q_k:\Spec(S/pS)\hookrightarrow\Spec(R/pR)$ introduced after the Theorem 2.2 and associated to the deformation package $\dbP(\scrB)$,  does not depend on $\scrB\in\scrB_E$.}

\medskip
\proof
Due to the property (ii), we have identities $n_0=n$, $N_0=N$, $U_0=U$, and $q_0=q$. Let $\bar S:=S/pS$. Let $\bar J:=J/pJ$; it is the maximal ideal of $\bar S$. Let $\scrB_1\in\scrB_E$. Let $\dbP(\scrB_1)=(n,N_1,n,N_1,U_1,U_1,u_{1\text{univ}},\Phi_{N_1},\nabla_1,\grH_{1S/pS},q_{1k},q_{1k},\grK_1,\bar\grK_1)$ be the deformation package of $\scrB_1$. We write $\nabla=\delta_0+\omega$ and $\nabla_1=\delta_0+\omega_1$, where $\omega,\omega_1\in E\otimes_{W(k)} \Omega_S^\wedge$ (cf. Proposition 2.6 (a)). As $\omega,\omega_1\in E\otimes_{W(k)} \Omega_S^\wedge$, due to Lemma 2.3 and the property (i) the reductions modulo $p$ of the two Kodaira--Spencer maps $\grK$ and $\grK_1$ have equal images. The Corollary is a direct consequence of this fact; for the sake of completeness, we present the details.   

We will check by induction on $m\in\dbN$ that there exists a unique automorphism $j_m:\Spec(\bar S/\bar J^m)\arrowsim \Spec(\bar S/\bar J^m)$ of $k$-schemes such that there exists an isomorphism $\theta_m:j_m^*(\grH\times_S \bar S/\bar J^m)\arrowsim\grH_1\times_S \bar S/\bar J^m$ which lifts the identity automorphism of $D$. As $\scrD_k=\Spf(R/pR)$ is the formal deformation space of $D$, this is equivalent to the fact that the closed subschemes $q_{1k,m}:\Spec(\bar S/\bar J^m)\hookrightarrow \Spec(R/pR)$ and $q_{k,m}:\Spec(\bar S/\bar J^m)\hookrightarrow \Spec(R/pR)$ defined by $q_{1k}$ and $q_k$ (respectively) are equal and we have $q_{1k,m}=q_{k,m}\circ j_m$; obviously this implies that if $j_m$ exists, then it is unique. The basis of the induction holds for $m=1$. The passage from $m$ to $m+1$ goes as follows. 

Let $\tilde j_{m+1}:\Spec(\bar S/\bar J^{m+1})\arrowsim \Spec(\bar S/\bar J^{m+1})$ be an arbitrary automorphism of $k$-schemes that lifts $j_m$. We endow the ideal $\bar J^m/\bar J^{m+1}$ of $\bar S/\bar J^{m+1}$ with the trivial divided power structure (see [Me, Ch. 5, Subsubsect. 2.3.5]; thus $(\bar J^m/\bar J^{m+1})^{[l]}=0$ for all $l\in\dbN\setminus\{1\}$). The evaluations at the trivial thickening $\Spec(\bar S/\bar J^{m+1})\arrowsim\Spec(\bar S/\bar J^{m+1})$ of the filtered $F$-crystals  of the two $p$-divisible groups $\tilde j_{m+1}^*(\grH\times_S \bar S/\bar J^{m+1})$ and $\grH_1\times_S \bar S/\bar J^{m+1}$ over $\bar S/\bar J^{m+1}$, can be identified with $(\bar M\otimes_k \bar S/\bar J^{m+1},\bar F^1\otimes_k \bar S/\bar J^{m+1},\bar\Phi_N(m+1),\bar\nabla(m+1))$ and respectively $(\bar M\otimes_k \bar S/\bar J^{m+1},F^1_{1,m+1},\bar\Phi_N(m+1),\bar\nabla(m+1))$, where $(\bar\Phi_N(m+1),\bar\nabla(m+1))$ is the reduction modulo $(p,J^{m+1})$ of $(\Phi_N,\nabla)$ and where $F^1_{1,m+1}$ is a direct summand of $\bar M\otimes_k \bar S/\bar J^{m+1}$ that lifts $\bar F^1\otimes_k \bar S/\bar J^m$. As we have $\omega,\omega_1\in E\otimes_{W(k)} \Omega_S^\wedge$, there exists an element $e_{m+1}\in E/pE\otimes_k \bar J^m/\bar J^{m+1}$ such that we have $F^1_{1,m+1}=(1_{\bar M\otimes_k \bar S/\bar J^{m+1}}+e_{m+1})(\bar F^1\otimes_k \bar S/\bar J^{m+1})$. The direct summand $F^1_{1,m+1}$ of $\bar M\otimes_k \bar S/\bar J^{m+1}$ depends only on the image $\bar e_{m+1}$ of $e_{m+1}$ in $\scrT\otimes_k \bar J^m/\bar J^{m+1}$. As $\omega,\omega_1\in E\otimes_{W(k)} \Omega_S^\wedge$, due to Lemma 2.3 and the property (i) the reductions modulo $(p,J^{m+1})$ of $\grK$ and $\grK_1$ have the same image $\nu(E)\otimes_k \bar S/\bar J^{m+1}$ which contains $\bar e_{m+1}$. Therefore we can replace $\tilde j_{m+1}$ by another automorphism $j_{m+1}:\Spec(\bar S/\bar J^{m+1})\arrowsim \Spec(\bar S/\bar J^{m+1})$ that lifts $j_m$, so that under this replacement the element $\bar e_{m+1}$ gets replaced by (i.e., becomes) the zero element. Thus we can assume that $e_{m+1}=0$; this implies that $F^1_{1,m+1}=\bar F^1\otimes_k \bar S/\bar J^{m+1}$. From the Grothendieck--Messing deformation theory applied to the trivial thickening $\Spec(\bar S/\bar J^m)\hookrightarrow \Spec(\bar S/\bar J^{m+1})$, we get that there exists a unique isomorphism $\theta_{m+1}:j_{m+1}^*(\grH\times_S \bar S/\bar J^{m+1})\arrowsim\grH_1\times_S \bar S/\bar J^{m+1}$ that lifts $\theta_m$. Thus $j_{m+1}$ exists and (as argued above) it is also unique. This ends the induction. 

As $\bar S$ is complete in the $\bar J$-adic topology, there exists a unique automorphism $j_{\infty}:\Spec(\bar S)\arrowsim\Spec(\bar S)$ that lifts each $j_m$ with $m\in\dbN$. The identities $q_{1k,m}=q_{k,m}\circ j_m$ with $m\in\dbN$ imply that $q_{1k}=q_k\circ j_{\infty}$. Therefore we conclude that the closed embedding $q_k:\Spec(S/pS)\hookrightarrow\Spec(R/pR)$ does not depend on the choice of an element $\scrB\in\scrB_E$.\endproof

\bigskip
\noindent
{\boldsectionfont 3. The proofs of the Basic Theorems A and B}
\bigskip 
 
In this Section we prove the two Basic Theorems A and B (see Subsections 3.3 and 3.5). Lemmas 3.1 and 3.2 are ingredients in the proofs. In particular, Lemma 3.2 pertains to the axioms 1.4 (i) and (ii). Proposition 3.6 checks that the two Basic Theorems A and B are functorial in a natural way. Example 3.4 and Corollary 3.7 show that each $\scrG_k^E$ as in the Basic Theorems A, is a formal Lie subgroup of a suitable $Ext^1$ formal group associated to $D$ and $E$. 

Let $n$ and the sextuple $(S,I,J,\Phi_S,\Omega_S^\wedge,d\Phi_S)$ be as in the beginning of Section 2. 

\medskip\smallskip\noindent
{\bf 3.1. Lemma.} {\it Let $E$ be a $W(k)$-submodule of $V_-$ such that the pair $(E,p\phi)$ is a Dieudonn\'e module over $k$. 

\medskip
{\bf (a)} Then the Newton polygon slopes of the $p$-divisible group $\dbD^{-1}(E)$ over $k$ belong to the interval $[0,1)$ (i.e., the $p$-divisible group $\dbD^{-1}(E)^{\text{t}}$ is connected).

\smallskip
{\bf (b)} Let $\scrB=(v_1,\ldots,v_n)\in E^n$ be such that $\{\nu(v_1),\ldots,\nu(v_n)\}$ generates a $k$-vector subspace of $\nu(E)$ of dimension $n$. Let $\dbP(\scrB)=(n,N,n,N,U,U,u_{\text{univ}},\Phi_N,\nabla,\grH_{S/pS},q_k,q_k,\grK,\bar\grK)$ be the deformation package of $\scrB$ (we note that $N\subseteq E$). Then $\grH_{S/pS}$ is a geometrically constant $p$-divisible group over $S/pS$.}

\medskip
\proof
As all Newton polygon slopes of $(V_-,\phi)$ belong to the interval $[-1,0)$, all Newton polygon slopes of $(E,p\phi)$ (i.e., of $\dbD^{-1}(E)$) belong to the interval $[0,1)$. Thus (a) holds. 

To check (b), let $K$, $h:\Spec(K)\to\Spec(S)$, and $u_h\in U(W(K))$ be as in the Fact 2.4. Let $n_h\in N\otimes_{W(k)} W(K)$ be such that $u_h=1_{M\otimes_{W(k)} W(K)}+n_h$. To prove the Lemma, we have to check that the pull back of $\grH_{S/pS}$ via the morphism $\Spec(K)\to\Spec(S/pS)$ defined by $h$ is isomorphic to $D_K$ i.e., the Dieudonn\'e modules $(M\otimes_{W(k)} W(K),u_h(\phi\otimes\sigma_K))$ and $(M\otimes_{W(k)} W(K),\phi\otimes\sigma_K)$ are isomorphic. The Verschiebung map of $(E,p\phi)$ is the restriction to $E$ of the $\sigma^{-1}$-linear automorphism $\phi^{-1}$ of $\End(M[{1\over p}])$. Thus for all $i\in\dbN$ we have $\phi^{-i}(E)\subseteq E$; in particular $n_{h,i}:=(\phi\otimes\sigma_K)^{-i}(n_h)\in E\otimes_{W(k)} W(K)$. As $E\subseteq V_-$, all Newton polygon slopes of $(E,\phi^{-1})$ are positive. Thus the sequence $(n_{h,i})_{i\in\dbN}$ of elements of $E\otimes_{W(k)} W(K)$ converges to $0$ in the $p$-adic topology. Let $u_{h,i}:=1_{M\otimes_{W(k)} W(K)}+n_{h,i}$. As $E\subseteq V_-$, each $n_{h,i}\in E\otimes_{W(k)} W(K)$ is a nilpotent element of $\End(M)\otimes_{W(k)} W(K)$. Thus we have $u_{h,i}\in \pmb{\GL}_M(W(K))$. The sequence $(u_{h,i})_{i\in\dbN}$ of elements of $\pmb{\GL}_M(W(K))$ converges to $1_{M\otimes_{W(k)} W(K)}$ in the $p$-adic topology. Thus the element $u_{h,\infty}:=\lim_{i\to\infty} u_{h,i}u_{h,i-1}\cdots u_{h,1}\in \pmb{\GL}_M(W(K))$
is well defined. For $g\in \pmb{\GL}_M(W(K))$, let $(\phi\otimes\sigma_K)(g):=\phi\otimes\sigma_K\circ g\circ (\phi\otimes\sigma_K)^{-1}\in \pmb{\GL}_M(B(K))$. Let $u_{h,0}:=u_h$. For $i\in\dbN$ we have $(\phi\otimes\sigma_K)(u_{h,i})=u_{h,i-1}\in \pmb{\GL}_M(W(K))$. Thus
$$u_{h,\infty}u_h (\phi\otimes\sigma_K)(u_{h,\infty})^{-1}=\lim_{i\to\infty} u_{h,i}\cdots u_{h,1}u_{h,0}(\phi\otimes\sigma_K)(u_{h,1})^{-1}\cdots (\phi\otimes\sigma_K)(u_{h,i})^{-1}$$
$$=\lim_{i\to\infty} u_{h,i}\cdots u_{h,0}u_{h,0}^{-1}\cdots u_{h,i-1}^{-1}=\lim_{i\to\infty} u_{h,i}=1_{M\otimes_{W(k)} W(K)}.$$
This implies that $u_{h,\infty}\in \pmb{\GL}_M(W(K))$ is an isomorphism between the Dieudonn\'e modules $(M\otimes_{W(k)} W(K),u_h(\phi\otimes\sigma_K))$ and $(M\otimes_{W(k)} W(K),\phi\otimes\sigma_K)$ over $K$.\endproof

\medskip\smallskip\noindent
{\bf 3.2. Lemma.} {\it  Let $W$ be a $B(k)$-vector subspace of $L_-$ which is stable under $\phi$.  Let $E:=\{x\in V_-\cap W|\phi^{-q}(x)\in V_-\cap W\;\forall q\in\dbN\}$. Let $E_0$ be a $W(k)$-submodule of $V_-\cap W$ such that the pair $(E_0,p\phi)$ is a Dieudonn\'e module. We have:

\medskip
{\bf (a)} Then $E$ is the largest $W(k)$-submodule of $V_-\cap W$ such that the pair $(E,p\phi)$ is a Dieudonn\'e module over $k$ and therefore we have $E_0\subseteq E$.

\smallskip
{\bf (b)} Let $c_E\in\dbN\cup\{0\}$ be the codimension of $(E,p\phi)$ (i.e., of $\dbD^{-1}(E)$). Then $\dim_k(\nu(E))=c_E$. Moreover, we have $E_0=E$ if and only if $\dim_k(\nu(E_0))=c_E$.

\smallskip
{\bf (c)} If $W^2=0$ (i.e., if $E^2=0$), then there exist:

\medskip\noindent
{\bf (c.i)} an element $t\in L_0$ which is fixed by $\phi$, which is a projector (i.e., $t^2=t$), and for which we have $[E,t]=E$;

\smallskip\noindent
{\bf (c.ii)} a short exact sequence $0\to (M_1,\phi)\to (M,\phi)\to (M_2,\phi)\to 0$ such that $E$ is contained in the direct summand $\Hom(M/M_1,M_1)$ of $\End(M)$.}

\medskip
\proof
The Verschiebung map for $p\phi$ is $\phi^{-1}$. As $p\phi(E)\subseteq V_-\cap W$ and for all $q\in\dbN$ we have $\phi^{-q}(p\phi(E))=p\phi^{-q+1}(E)\subseteq V_-\cap W$, we also have $p\phi(E)\subseteq E$. Obviously $\phi^{-1}(E)\subseteq E$. Therefore the pair $(E,p\phi)$ is a Dieudonn\'e module over $k$. As we have $\phi^{-q}(E_0)\subseteq E_0\subseteq V_-\cap W$ for all $q\in\dbN\cup\{0\}$, we have $E_0\subseteq E$. From this (a) follows. 

We prove (b). Let $x\in E\setminus pE$ and $y\in E\setminus pE$ be such that $p\phi(x)\in E\setminus pE$ and $p\phi(y)\in pE$. We check that $\nu(x)\in\scrT$ is a non-zero element and that $\nu(y)=0$. We show that the assumption that $\nu(x)=0$ leads to a contradiction. As $\nu(x)=0$ we have $\phi(x)\in\End(M)\cap W=V_-\cap W$, cf. property 1.3.1 (*). For all $q\in\dbN\cup\{0\}$ we have $\phi^{-q}(\phi(x))=\phi^{-q+1}(x)\subseteq V_-\cap W$. Thus $\phi(x)\in E$. Therefore $p\phi(x)\in pE$ and this contradicts the assumption that $p\phi(x)\in E\setminus pE$. Thus $\nu(x)\neq 0$. As $p\phi(y)\in pE$, we have $\phi(y)\in E\subseteq V_-\cap W$. But if $\nu(y)\neq 0$, then $\phi(y)\notin\End(M)$ (cf. property 1.3.1 (*)). Therefore $\nu(y)=0$. Let $\bar E:=E/pE$ and let $\bar F^1_E$ be the kernel of the reduction modulo $p$ of $p\phi:E\to E$. The properties of $\nu(x)$ and $\nu(y)$ imply that the map $\bar E\to\scrT$ induced by $\nu$ has $\bar F^1_E$ as its kernel. Thus the image of $E$ in $\scrT$ can be naturally identified with $\bar E/\bar F_E^1$ and therefore it is a $k$-vector space of dimension $c_E$. 

To end the proof of (b) it suffices to show that the equality $\dim_k(\nu(E_0))=c_E$ implies that $E=E_0$. Let $\eps:\dbD^{-1}(E)\to \dbD^{-1}(E_0)$ be the homomorphism of $p$-divisible groups over $k$ defined naturally by the monomorphism $(E_0,p\phi)\hookrightarrow (E,p\phi)$ of Dieudonn\'e modules over $k$. Let $\iota:\dbD^{-1}(E_0)^{\text{t}}\to \dbD^{-1}(E)^{\text{t}}$ be the Cartier dual of $\eps$. Let $\bar E_0:=E_0/pE_0$. As $\dim_k(\nu(E_0))=c_E$, we have $\nu(E_0)=\nu(E)$. Thus the natural $k$-linear map $\bar E_0\to \bar E/\bar F^1_E$ is onto. Therefore, if $\bar F^1_{E_0}$ is the kernel of the reduction modulo $p$ of $p\phi:E_0\to E_0$, then we have a $k$-linear epimorphism $\bar E_0/\bar F^1_{E_0}\twoheadrightarrow \bar E/\bar F^1_E$ which is the tangent map $d\iota$ of $\iota$. From this and the fact that $\dbD^{-1}(E)^{\text{t}}$ is connected (see Lemma 3.1 (a)), we get that $\iota$ is an epimorphism. From this and the fact that the height of $\dbD^{-1}(E)^{\text{t}}$ is at least equal to the height of $\dbD^{-1}(E_0)^{\text{t}}$, we get that $\iota$ is an isogeny. Thus $d\iota$ is a $k$-linear epimorphism between $k$-vector spaces of the same dimension and therefore it is a $k$-linear isomorphism. This and the fact that $\iota:\dbD^{-1}(E_0)^{\text{t}}\to \dbD^{-1}(E)^{\text{t}}$ is an isogeny between connected $p$-divisible groups, implies that $\iota$ is an isomorphism. Thus $E_0=E$. Therefore (b) holds. 

We prove (c). Let $W_0$ be the $B(k)$-vector subspace of $M[{1\over p}]$ generated by all $y(M[{1\over p}])$ with $y\in W$. If $\beta_W$ is a $B(k)$-basis for $W$, then we have $W_0=\sum_{y\in\beta_W} y(M[{1\over p}])$. As $\phi(W)=W$, we have $\phi(W_0)=W_0$. We consider a direct sum decomposition $(M[{1\over p}],\phi)=(W_0,\phi)\oplus (W_0^\perp,\phi)$ into $F$-isocrystals. Let $t\in L_0$ be the projection of $M[{1\over p}]$ on $W_0^\perp$ along $W_0$. For each $x\in W$, we have $x(M[{1\over p}])\subseteq W_0$ (cf. the very definition of $W_0$). Moreover $x(W_0)=\sum_{y\in\beta_W} xy(M[{1\over p}])=0$ (as each $xy\in W^2=0$). From the last two sentences we get that $tx=0$, $xt=x$, and $[x,t]=x$. Thus $[E,t]=E$. Therefore (c.i) holds. If $M_1:=M\cap W_0$ and $M_2:=M/M_1$, then we get a short exact sequence $0\to (M_1,\phi)\to (M,\phi)\to (M_2,\phi)\to 0$ such that $W$ is contained in $\Hom(M/M_1,M_1)[{1\over p}]$. This implies that $E\subseteq \Hom(M/M_1,M_1)$ and thus (c.ii) also holds.\endproof

\medskip\smallskip\noindent
{\bf 3.3. Proof of the Basic Theorem A.} The proof of Basic Theorem A is divided into five parts (see the below Subsubsections 3.3.1 to 3.3.5). 

\medskip
$\bullet$ The first part introduces a deformation package $\dbP(\scrB)$ that produces a closed embedding $q_k:\Spec(S/pS)\hookrightarrow \Spec(R/pR)$ which will define $\scrG_k^E$. 

$\bullet$ The second part uses Lemma 3.2 (c) to produce an enlargement $\tilde E$ of $E$ which is a Lie subalgebra of $\End(M)[{1\over p}]$ equipped with a short exact sequence $0\to E\to\tilde E\to W(k)\to 0$. The role of $\tilde E$ is to produce from $\dbP(\scrB)$ a suitable short exact sequence of $F$-crystals over $S/pS$ (for this we will use the connection results of Section 2). 

$\bullet$ The third part shows that the short exact sequence of $F$-crystals over $S/pS$ is (isomorphic to) the one associated to the universal extension of $\dbD^{-1}(E)$ by $(\pmb{\mu}_{p^{\infty}})_k$. This endows $\Spf(S/pS)$ with a natural structure of a formal Lie group.

$\bullet$ The fourth part checks that if $\scrG_k^E\hookrightarrow \scrD_k$ is indeed defined by the closed embedding $q_k:\Spec(S/pS)\hookrightarrow \Spec(R/pR)$, then the non-canonical parts of Basic Theorem A hold. Lemma 3.1 (b) is used in this part. 

$\bullet$ The fifth part checks that $\scrG_k^E$ and its formal Lie group structure are indeed canonical i.e., do not depend on the choices of $\dbP(\scrB)$ and $\tilde E$ made in the first two parts.  

\medskip\noindent
{\bf 3.3.1. Explicit deformation package.}  We take $n\in\dbN\cup\{0\}$ to be the codimension of $\dbD^{-1}(E)$. Let $\scrB=(v_1,\ldots,v_n)$ be an $n$-tuple of elements of $E$ such that $\{\nu(v_1),\ldots,\nu(v_n)\}$ is a $k$-basis for $\nu(E)$, cf. Lemma 3.2 (b). We consider the deformation package $\break\dbP(\scrB)=(n,N,n,N,U,U,u_{\text{univ}},\Phi_N,\nabla,\grH_{S/pS},q_k,q_k,\grK,\bar\grK)$ of $\scrB$. Let $p\tilde\Phi_N:\End(M)\otimes_{W(k)} S\hookrightarrow \End(M)\otimes_{W(k)} S$ be the $\Phi_S$-linear monomorphism induced by $p\Phi_N$. Let $l_U:U\hookrightarrow \tilde U$ be the closed embedding of the proof of Fact 2.5. We have $p\tilde\Phi_N=\tilde u_{\text{univ}}(p\phi\otimes\Phi_S)$, where $\tilde u_{\text{univ}}:=1_{\tilde E\otimes_{W(k)} S}+\sum_{i=1}^n dl_U(v_i)\otimes x_i\in \tilde U(S)\leqslant \pmb{\GL}_{\tilde E}(S)$ (cf. Formulas (1) and (2)). The morphism $q_k:\Spec(S/pS)\to\Spec(R/pR)$ is a closed embedding, cf. Lemma 2.3. 

\medskip\noindent
{\bf 3.3.2. Lie elements and connections.} Let $t\in L_0$ be an element that is fixed by $\phi$ and such that we have $[E,t]=E$, cf. Lemma 3.2 (c.i). We refer to $t$ as a {\it Lie element} for $E$. Let $\tilde E:=E\oplus W(k)t$. As $E^2=0$ and $[E,t]=E$, $\tilde E$ is a Lie subalgebra of $\End(M[{1\over p}])$. As $N=N_0\subseteq E$, we have $N^2=0$ and $N$ is a direct summand of $\End(M)$. Thus $U$ is a smooth, closed subgroup scheme of $\pmb{\GL}_M$ whose Lie algebra is $N$ (cf. properties (i) and (ii) before Theorem 2.2). It is convenient to identify $\tilde E/E=W(k)t$, whenever we have short exact sequences that are natural tensorizations of the short exact sequence $0\to E\to\tilde E\to \tilde E/E\to 0$ (of either Lie algebras over $W(k)$ or $W(k)$-modules). 

We have $\nabla=\delta_0+\omega$ with $\omega\in E\otimes_{W(k)} \Omega_S^\wedge$, cf. Proposition 2.6 (a). The connection $\nabla$ induces a connection $\tilde\nabla=\delta_0+\tilde\omega$ on $\tilde E\otimes_{W(k)}  S$. As $\tilde\omega\in \Hom(\tilde E/E,E)\otimes_{W(k)} \Omega_S^\wedge$ (cf. Proposition 2.6 (b)), the connection $\tilde\nabla$ induces the flat connection $\delta_0$ on both $E\otimes_{W(k)} S$ and $\tilde E/E\otimes_{W(k)} S=St$. Thus we have a short exact sequence of $F$-crystals over $S/pS$
$$0\to (E\otimes_{W(k)}  S,p\phi\otimes\Phi_S,\delta_0)\to (\tilde E\otimes_{W(k)}  S,p\tilde\Phi_N,\tilde\nabla)\to (St,p\phi\otimes\Phi_S,\delta_0)\to 0.\leqno (6)$$
As the left $U$-modules $E$ and $\tilde E/E$ are trivial, for the first and the last terms of (6) we do not have to insert the universal element $u_{\text{univ}}\in U(S)\leqslant \pmb{\GL}_M(S)$. Let $\grL_{S/pS}$ be the unique $p$-divisible group over $\Spec(S/pS)$ whose $F$-crystal over $S/pS$ is $(\tilde E\otimes_{W(k)}  S,p\tilde\Phi_N,\tilde\nabla)$ (its existence and uniqueness is argued in the same way we argued the existence and the uniqueness of $\grH_{S/pS}$ in the paragraph after Theorem 2.2, the role of $\scrB=\{v_1,\ldots,v_n\}$ being replaced by $\{dl_U(v_1),\ldots,dl_U(v_n)\}$). As $\{x_1,\ldots,x_n\}$ is a $p$-basis for the regular ring $S/pS$, the crystalline Dieudonn\'e functor $\dbD$ is fully faithful on $p$-divisible groups over $\Spec(S/pS)$ (cf. [BM, Thm. 4.1.1]). Thus to (6) corresponds a complex
$$0\to (\pmb{\mu}_{p^{\infty}})_{S/pS}\to \grL_{S/pS}\to \dbD^{-1}(E)\times_k S/pS\to 0\leqno (7)$$
of $p$-divisible groups over $\Spec(S/pS)$. The complex (7) is a short exact sequence as this is so modulo the maximal ideal of $S/pS$. The Cartier dual of (7) is a short exact sequence
$$0\to \dbD^{-1}(E)^{\text{t}}\times_k {S/pS}\to \grL^{\text{t}}_{S/pS}\to (\dbQ_p/\dbZ_p)_{S/pS}\to 0\leqno (8)$$
of $p$-divisible groups over $\Spec(S/pS)$.
\medskip\noindent
{\bf 3.3.3. Lemma.} {\it The short exact sequence (7) is the universal extension of $\dbD^{-1}(E)$ by $(\pmb{\mu}_{p^{\infty}})_k$. In other words, if $Y$ is a local, artinian $k$-scheme of residue field $k$ and if we have a short exact sequence  $0\to (\pmb{\mu}_{p^{\infty}})_Y\to\grE\to \dbD^{-1}(E)\times_k Y\to 0$ of $p$-divisible groups over $Y$ that lifts the split short exact sequence $0\to (\pmb{\mu}_{p^{\infty}})_k\to \dbD^{-1}(E)\oplus (\pmb{\mu}_{p^{\infty}})_k\to \dbD^{-1}(E)\to 0$, then there exists a unique morphism $Y\to\Spec(S/pS)$ of $k$-schemes such that the pull back of (7) via it is (up to a unique isomorphism that lifts the identity automorphism of $\dbD^{-1}(E)\oplus (\pmb{\mu}_{p^{\infty}})_k$) the short exact sequence $0\to (\pmb{\mu}_{p^{\infty}})_Y\to\grE\to \dbD^{-1}(E)\times_k Y\to 0$.}
\medskip
\proof
As the map $N/pN\to\scrT$ induced by $\nu$ is injective, from the property 1.3.1 (*) we get that $p\phi(N)$ is a direct summand of $\End(M)$ and thus also of $E$. As the rank $n$ of $N$ is the codimension of $\dbD^{-1}(E)$, there exists a direct sum decomposition $E=N\oplus F^1_E$ into $W(k)$-modules such that the triple $(E,F^1_E,p\phi)$ is a filtered Dieudonn\'e module. It is well known that a universal extension of $\dbD^{-1}(E)$ by $(\pmb{\mu}_{p^{\infty}})_k$ exists and it is associated naturally to a smooth formal subscheme of the formal deformation space of $\dbD^{-1}(E)\oplus (\pmb{\mu}_{p^{\infty}})_k$ over $\Spf(k)$. The tangent space of the universal extension of $\dbD^{-1}(E)$ by $(\pmb{\mu}_{p^{\infty}})_k$ (i.e., of $(W(k)t,p\phi)$ by $(E,p\phi)$), is the $k$-vector subspace 
$$\scrT^{\text{ue}}_0:=\Hom(\tilde E/E,E/F^1_E)\otimes_{W(k)} k=\Hom((F^1_E+W(k)t)/F^1_E,E/F^1_E)\otimes_{W(k)} k$$ of the tangent space 
$$\scrT_0:=\Hom(F^1_E+W(k)t,E/F^1_E)\otimes_{W(k)} k=\End(\tilde E/p\tilde E)/F^0(\End(\tilde E/p\tilde E))$$ 
of the formal deformation space of the direct sum $p$-divisible group $\dbD^{-1}(E)\oplus (\pmb{\mu}_{p^{\infty}})_k$; here $\break F^0(\End(\tilde E/p\tilde E)):=\{x\in\End(\tilde E/p\tilde E)|x(F^1_E/pF^1_E+kt)\subseteq F^1_E/pF^1_E+kt\}$. Therefore $\dim_k(\scrT^{\text{ue}}_0)=n$. Thus (7) is the universal extension of $\dbD^{-1}(E)$ by $(\pmb{\mu}_{p^{\infty}})_k$ if and only if the tangent space $\scrT_0^{\text{b}}$ of (7) (i.e., the image of the tensorization with the maximal ideal $I$ of $S$ of the Kodaira--Spencer map of $\tilde\nabla$)  has dimension $n$.

The image of the differential map $dl_U:\Lie(U)\to \Lie(\tilde U)$ defined by Formula (2) is $\im(dl_U)=\Hom(\tilde E/E,[N,t])$. As $\tilde\nabla$ is induced by $\nabla$ and as $N=\Lie(U)$, from Lemma 2.3 we get that the tangent space $\scrT_0^{\text{b}}$ is the projection into $\scrT_0$ of $\im(dl_U)$. As $\phi(t)=t$ and $[E,t]=E$, the association $x\to [x,t]$ defines an automorphism $(E,p\phi)\arrowsim (E,p\phi)$. Therefore we have another direct sum decomposition $E=[N,t]\oplus F^1_E$ into $W(k)$-modules. As the $k$-vector subspace $[N,t]\otimes_{W(k)} k$ of $E\otimes_{W(k)} k$  maps isomorphically onto $(E/F^1_E)\otimes_{W(k)} k$, the image $\scrT_0^{\text{b}}$ of $\im(dl_U)$ in $\scrT_0$ is exactly $\scrT^{\text{ue}}_0$. Thus $\scrT_0^{\text{b}}=\scrT^{\text{ue}}_0$. Therefore $\dim_k(\scrT_0^{\text{b}})=n$. Thus the short exact sequence (7) is the universal extension of $\dbD^{-1}(E)$ by $(\pmb{\mu}_{p^{\infty}})_k$.\endproof

\medskip\noindent
{\bf 3.3.4. Formal Lie group structure.} As (8) is the Cartier dual of (7), from Lemma 3.3.3 we get that (8) is the universal extension of $(\dbQ_p/\dbZ_p)_k$ by $\dbD^{-1}(E)^{\text{t}}$. Thus the formal subscheme $\Spf(S/pS)$ of $\scrD_k=\Spf(R/pR)$ has a natural structure of a formal Lie group over $\Spf(k)$ isomorphic to the formal Lie group $Ext^1((\dbQ_p/\dbZ_p)_k,\dbD^{-1}(E)^{\text{t}})$ over $\Spf(k)$. As $Hom((\dbZ_p)_k,\dbD^{-1}(E)^{\text{t}})\arrowsim Ext^1((\dbQ_p/\dbZ_p)_k,\dbD^{-1}(E)^{\text{t}})$ (cf. third paragraph of Section 1), the formal subscheme $\Spf(S/pS)$ of $\scrD_k=\Spf(R/pR)$ has a natural structure of a formal Lie group over $\Spf(k)$ isomorphic to the formal Lie group over $\Spf(k)$ associated to the $p$-divisible group $\dbD^{-1}(E)^{\text{t}}$. The addition law of $\Spf(S/pS)$ is defined naturally by adding the two short exact sequences over $\Spec(\widehat{S/pS\otimes_k S/pS})$ that are the pulls back of (8) via the two projections $s_1,s_2:\Spec(\widehat{S/pS\otimes_k S/pS})\to\Spec(S/pS)$ (here $\widehat{S/pS\otimes_k S/pS}$ is the completion of the local ring $S/pS\otimes_k S/pS$). As $\scrG_k^E$ we take the formal subscheme $\Spf(S/pS)$ of $\scrD_k=\Spf(R/pR)$ defined by the closed embedding $q_k:\Spec(S/pS)\hookrightarrow \Spec(R/pR)$, together with its formal Lie group structure we have introduced. From the last three sentences we get that the non-canonical part of the property 1.5 (i) holds. From Lemma 2.3 and the identities $N_0=N$ and $\nu(E)=\nu(N)$, we get that the property 1.5 (ii) holds. The property 1.5 (iii) also holds, cf. Lemma 3.1 (b).

\medskip\noindent
{\bf 3.3.5. Canonical properties.}
Based on the previous Subsubsection, to end the proof of the Basic Theorem A (i.e., to check the canonical parts of it), we only have to show that the following three properties hold: 

\medskip
{\bf (i)} the formal Lie group structure on $\scrG_k^E$ does not depend on the choice of the Lie element $t\in L_0$ for $E$;

\smallskip
{\bf (ii)} the formal subscheme $\scrG_k^E$ of $\scrD_k$ does not depend on the choice of $\scrB\in E^{\times n}$;

\smallskip
\smallskip
{\bf (iii)} the formal Lie group structure on $\scrG_k^E$ does not depend on the choice of $\scrB\in E^{\times n}$.

\medskip\noindent
{\bf On (i).} To check (i), we consider another element $t_1\in L_0$ which is fixed by $\phi$ and for which we have $[E,t_1]=E$. Let $\tilde E_1:=E\oplus W(k)t_1$.  We have a Lie isomorphism 
$$j:\tilde E_1\arrowsim \tilde E$$ 
that is compatible with $p\phi$ and that obeys the following two rules: (i.a) $j(t_1)=t$, and (i.b) $j([x,t_1])=[x,t]$ for all $x\in E$. The Lie isomorphism $j$ is also an isomorphism of left $U$-modules. Thus $p\tilde\Phi_N\circ j\otimes 1_S=j\otimes 1_S\circ p\tilde\Phi_N:\tilde E_1\otimes_{W(k)} S\to \tilde E\otimes_{W(k)} S$. Let $\tilde\nabla_1$ be the connection on $\tilde E_1\otimes_{W(k)} S$ that is the analogue of $\tilde\nabla$ but obtained working with $t_1$ instead of $t$. As $\tilde\nabla_1$ is uniquely determined by $p\tilde\Phi_N$ restricted to $\tilde E_1\otimes_{W(k)} S$ (this is only a variant of Theorem 1.2 (a)) and as $p\tilde\Phi_N\circ j\otimes 1_S=j\otimes 1_S\circ p\tilde\Phi_N$, $j\otimes 1_S$ is parallel with respect to $\tilde\nabla_1$ and $\tilde\nabla$. In other words, $j\otimes 1_S$ induces an isomorphism between the short exact sequence (6) of $F$-crystals over $S/pS$ and the analogous short exact sequence of $F$-crystals over $S/pS$ one gets by working with $t_1$ instead of $t$. Therefore the isomorphism class of the short exact sequence (6) does not depend on the choice of $t$. From this and [BM, Thm. 4.1.1] we get that the isomorphism class of (7) and thus also of (8) does not depend on $t$. Therefore the formal Lie group structure on $\scrG^E_k$ does not depend on the choice of $t$. Thus (i) holds. 

\medskip\noindent
{\bf On (ii)}. The fact that (ii) holds is implied by Corollary 2.7. We have two extra properties: 

\medskip
{\bf (iv)} the $S$-submodule $E\otimes_{W(k)} S$ of $\End(M)\otimes_{W(k)} S$ is the unique $S$-submodule which modulo $J$ is $E$ and which defines a constant $F$-crystal $(E\otimes_{W(k)} S,p\phi\otimes\Phi_S,\delta_0)$ over $S/pS$ (in fact $L_-\otimes_{B(k)} S[{1\over p}]$ is the maximal $S[{1\over p}]$-submodule of $\End(M[{1\over p}])\otimes_{B(k)} S[{1\over p}]$ with the property that all slopes of $(L_-\otimes_{B(k)} S[{1\over p}],p\tilde\Phi_N)$ belong to the interval $[0,1)$, cf. the fact that the $p$-divisible group $\grH_{S/pS}=q_k^*(\grD)$ is geometrically constant); 

\smallskip
{\bf (v)} the element $t\in (\End(M[{1\over p}])/E[{1\over p}])\otimes_{B(k)} S[{1\over p}]$ is uniquely determined by the facts that it is fixed by $\tilde\Phi_N$, it is annihilated by $\tilde\nabla$, and it lifts the element $t$ of $(\End(M[{1\over p}])/E[{1\over p}],\phi)$ (this is so, as the inclusion $\Phi_S(J)\subseteq J^p$ implies that there exist no element in $\break (\End(M[{1\over p}])/E[{1\over p}])\otimes_{B(k)} J[{1\over p}]$ that is fixed by $\tilde\Phi_N$).

\medskip\noindent
{\bf On (iii).} Based on (iv) and (v) we get that the $S$-submodule $\tilde E\otimes_{W(k)} S$ of $\End(M)\otimes_{W(k)} S$ is intrinsically associated to the pair $(t,E)$ and to the $F$-crystal $(M\otimes_{W(k)} S,\Phi_N,\nabla)$ over $S/pS$. Due to the property (ii), both the closed embedding $q_k:\Spec(S/pS)\hookrightarrow\Spec(R/pR)$ and the last $F$-crystal over $S/pS$ are intrinsically associated to $E$. As the formal Lie group structure on the formal subscheme $\scrG^E_k$ of $\scrD_k$ was intrinsically defined in terms of $\tilde E\otimes_{W(k)} S$, we conclude (based on the last three sentences) that (iii) holds.\endproof

\medskip\smallskip\noindent
{\bf 3.4. Example.} Suppose that we have a direct sum decomposition $(M,\phi)=(M_1,\phi)\oplus (M_2,\phi)$ such that the greatest Newton polygon slope of $(M_2,\phi)$ is greater than the smallest Newton polygon slope of $(M_1,\phi)$. Let $D=D_1\times_k D_2$ be the corresponding product decomposition of $p$-divisible groups over $k$. Let $E$ be the largest $W(k)$-submodule of $\Hom(M_2,M_1)$ such that the pair $(E,p\phi)$ is a Dieudonn\'e module over $k$; our assumption on Newton polygon slopes implies that $E\neq 0$. Let $t$ be the projection of $M$ on $M_2$ along $M_1$. We have $\phi(t)=t$ and $[E,t]=E$. Let $\scrE_k$ be the formal subscheme of $\scrD_k$ that parametrizes extensions of $D_1$ by $D_2$. We can identify $\scrE_k=Ext^1(D_1,D_2)$; thus $\scrE_k$ has a canonical structure of a formal Lie group over $\Spf(k)$. As we have $\nabla=\delta_0+\omega$ with $\omega\in E\otimes_{W(k)} \Omega_S^\wedge$ (see Subsubsections 3.3.1 and 3.3.2), the connection $\nabla$ annihilates both $M_1\otimes 1$ and $M_2\otimes 1$ and it leaves invariant the $S$-submodule $\Hom(M_2,M_1)\otimes_{W(k)} S$ of $\End(M)\otimes_{W(k)} S$. The group scheme $U$ fixes both $M_1$ and $M_2$. From the last two sentences we get that there exist short exact sequences 
$$0\to (M_1\otimes_{W(k)} S,\phi\otimes\Phi_S,\delta_0)\to (M\otimes_{W(k)} S,\Phi_N,\nabla)\to (M_2\otimes_{W(k)} S,\phi\otimes\Phi_S,\delta_0)\to 0\;\;\text{and}\leqno (9)$$
$$0\to (\Hom(M_2,M_1)\otimes_{W(k)} S,p\phi\otimes\Phi_S,\delta_0)\to ((\Hom(M_2,M_1)\oplus W(k)t)\otimes_{W(k)} S,p\tilde\Phi_N,\nabla)$$
$$\to (St,p\phi\otimes\Phi_S,\delta_0)\to 0\leqno (10)$$
of $F$-crystals over $S/pS$. To the short exact sequence (9) corresponds a short exact sequence $0\to D_2\times_k S/pS\to \grH_{S/pS}\to D_1\times_k S/pS\to 0$ of $p$-divisible groups over $\Spec(S/pS)$, cf. [BM, Thm. 4.1.1]. Thus $\scrG^E_k$ is a smooth formal subscheme of $\scrE_k$. 

Let $\scrF_k$  be the formal group over $k$ that parametrizes isomorphism classes of extensions of the $F$-crystal $(M_1,\phi)$ over $k$ by the $F$-crystal $(M_2,\phi)$ over $k$. Let $\scrW_k$  be the formal group over $k$ that parametrizes isomorphism classes of extensions of the $F$-crystal $(W(k)t,p\phi)$ over $k$ by the $F$-crystal $(\Hom(M_2,M_1),p\phi)$ over $k$. 

Here are three key functors.

\medskip\noindent
$\bullet$ The crystalline Dieudonn\'e functor $\dbD$ defines a homomorphism $\dbD_k:\scrE_k\to\scrF_k$ of formal groups over $\Spf(k)$. It is a monomorphism, cf. Grothendieck--Messing deformation theory.

\smallskip\noindent
$\bullet$ The short exact sequence (10) is the push forward of (6) via the monomorphism 
$$(E\otimes_{W(k)} S,p\phi\otimes\Phi_S,\delta_0)\hookrightarrow (\Hom(M_2,M_1)\otimes_{W(k)} S,p\phi\otimes\Phi_S,\delta_0)$$ 
of (constant) $F$-crystals over $S/pS$. This property is a direct consequence of the fact that we have $\nabla=\delta_0+\omega$ with $\omega\in E\otimes_{W(k)} \Omega_S^\wedge$ (and it holds without assuming that $E[{1\over p}]=\Hom(M_2,M_1)[{1\over p}]$).  Thus we have a natural homomorphism $\dbE_k:\scrG^E_k\to\scrW_k$. 

\smallskip\noindent
$\bullet$ We have a natural monomorphism $\dbT_k:\scrF_k\hookrightarrow\scrW_k$ defined as follows.  Not to introduce extra notations, we will describe here only the restriction of $\dbT_k$ to $\im(\dbE_k)$ (the general description of $\dbT_k$ is the same). The short exact sequence (10) is obtained from (9) via the following two operations (performed in the order listed):

\medskip
{\bf (i)} the tensorization of (9) with the $F$-crystal  $(\Hom(M_2,W(k))\otimes_{W(k)} S,p\phi\otimes\Phi_S,\delta_0)$ over $S/pS$ of the $p$-divisible group $D_2^{\text{t}}\times_k S/pS$;

\smallskip
{\bf (ii)} a pull back via the canonical embedding of $(W(k)t\otimes_{W(k)} S,p\phi\otimes\Phi_S,\delta_0)$ into the tensor product $(M_2\otimes_{W(k)} S,\phi\otimes\Phi_S,\delta_0)\otimes (\Hom(M_2,W(k))\otimes_{W(k)} S,p\phi\otimes\Phi_S,\delta_0)$ (the projection $t$ of $M$ on $M_2$ along $M_1$ maps canonically to $1_{M_2}\in M_2\otimes_{W(k)}\Hom(M_2,W(k))$).

\medskip\noindent
As the operations (i) and (ii) are compatible with additions of extensions, the restriction of $\dbT_k$ to $\im(\dbE_k)$ is a homomorphism of formal groups over $\Spf(k)$. 

\medskip
The composite of the formal closed embedding $\scrG_k^E\hookrightarrow\scrE_k$ with the monomorphism $\dbT_k\circ \dbD_k:\scrE_k\hookrightarrow\scrW_k$ is the homomorphism $\dbE_k:\scrG^E_k\to\scrW_k$. Thus $\scrG^E_k$ is a formal Lie subgroup of $\scrE_k$ as well as of $\scrW_k$. If both $D_1$ and $D_2$ are isoclinic $p$-divisible groups, then this result is also claimed in [Ch2].

\medskip\smallskip\noindent
{\bf 3.5. Proof of the Basic Theorem B.} We will check one by one that properties 1.6 (i) to (iv) hold.

\medskip\noindent
{\bf On 1.6 (i).} If $p>2$, then the existence and the uniqueness of $\dbD^{-1}(E)_{W(k)}$ is implied by the Grothendieck--Messing deformation theory. If $p=2$, then the existence and the uniqueness of $\dbD^{-1}(E)_{W(k)}$ is (for instance) a consequence of [Fo, Ch. IV, Prop. 1.6]. More precisely, we know that $\dbD^{-1}(E)^{\text{t}}$ is connected (cf. Lemma 3.1 (a)) and loc. cit. asserts that the category of $p$-divisible groups over $W(k)$ whose Cartier duals are connected is antiequivalent to the category of filtered Dieudonn\'e modules over $k$ that do not have slope $1$. Strictly speaking, loc. cit. is stated in terms of Honda triples (of the form $(M,\phi({1\over p}F^1),\phi_1)$) and not in terms of filtered Dieudonn\'e modules (of the form $(M,F^1,\phi_1)$). Thus regardless of the parity of $p$, the property 1.6 (i) holds.

\medskip\noindent
{\bf On 1.6 (ii).} Let the direct sum decomposition $E=F^0(E)\oplus F^{-1}(E)$ be as in the axiom 1.4 (iii). We have a canonical choice for $N$; it  is $N:=F^{-1}(E)$. Due to the axioms 1.4 (iii) and (iv), the short exact sequence (6) extends to a short exact sequence
$$0\to (E\otimes_{W(k)}  S,F^0(E)\otimes_{W(k)} S,p\phi\otimes\Phi_S,\delta_0)\to (\tilde E\otimes_{W(k)}  S,F^0(E)\otimes_{W(k)} S\oplus St,p\tilde\Phi_N,\tilde\nabla)$$$$\to (St,St,p\phi\otimes\Phi_S,\delta_0)\to 0\leqno (11)$$  
of filtered $F$-crystals over $S/pS$. For all $m\in\dbN$, the ideal $pJ$ of $S$ has a nilpotent divided power structure modulo $pJ^m$. From this and the Grothendieck--Messing deformation theory we get that to (11) corresponds a short exact sequence
$$0\to (\pmb{\mu}_{p^{\infty}})_S\to \grL\to \dbD^{-1}(E)_{W(k)}\times_{W(k)} S\to 0\leqno (12)$$
of $p$-divisible groups over $\Spec(S)$ that splits modulo $J$ and that lifts (7) (it is over $\Spec(S)$ and not only over $\Spf(S)$, cf. [Me, Ch. II, Lemma 4.16]). The Cartier dual of (12) is a short exact sequence of $p$-divisible groups over $\Spec(S)$ that is of the form
$$0\to \dbD^{-1}(E)_{W(k)}^{\text{t}}\times_{W(k)} S\to \grL^{\text{t}}\to (\dbQ_p/\dbZ_p)_S\to 0,\leqno (13)$$
that splits modulo $J$, and that lifts (8). Due to the Lemma 3.3.3, the short exact sequences (12) and (13) are universal extensions. Thus (13) is the universal extension of $(\dbQ_p/\dbZ_p)_{W(k)}$ by $\dbD^{-1}(E)_{W(k)}^{\text{t}}$. Therefore the formal subscheme $\Spf(S)$ of $\scrD=\Spf(R)$ has a  natural structure of a formal Lie group over $\Spf(W(k))$ isomorphic to the formal Lie group $Ext^1((\dbQ_p/\dbZ_p)_{W(k)},\dbD^{-1}(E)_{W(k)}^{\text{t}})$ over $\Spf(W(k))$. As we have an isomorphism of formal Lie groups $Hom((\dbZ_p)_{W(k)},\dbD^{-1}(E)_{W(k)}^{\text{t}})\arrowsim Ext^1((\dbQ_p/\dbZ_p)_{W(k)},\dbD^{-1}(E)_{W(k)}^{\text{t}})$ (see the third paragraph of Section 1), the formal subscheme $\Spf(S)$ of $\scrD=\Spf(R)$ has a natural structure of a formal Lie group over $\Spf(W(k))$ isomorphic to the formal Lie group over $\Spf(W(k))$ associated to the $p$-divisible group $\dbD^{-1}(E)_{W(k)}^{\text{t}}$. As $\scrG^E$ we take the formal subscheme of $\scrD=\Spf(R)$ defined by the closed embedding $q:\Spec(S)\hookrightarrow \Spec(R)$, together with its formal Lie group structure we have introduced. Thus the property 1.6 (ii) holds. Here $q$ is the morphism introduced after Theorem 2.2 and it is a closed embedding as its special fibre $q_k:\Spec(S/pS)\to\Spec(R/pR)$ is so.  

\medskip\noindent
{\bf On 1.6 (iii).} We recall that to the Hodge filtration $F^1$ of $D_{W(k)}$ corresponds a direct sum decomposition $(\tilde E,F^0(E)\oplus W(k)t,p\phi)=(E,F^0(E),p\phi)\oplus (W(k)t,W(k)t,p\phi)$ of filtered Dieudonn\'e modules over $k$. Thus as (13) splits  modulo $J$, the origin of $\scrG^E$ corresponds to $D_{W(k)}$ i.e., to the morphism $z_0:\Spec(W(k))\to\Spec(S)$ defined by the $W(k)$-epimorphism $w_0:S\twoheadrightarrow W(k)$ whose kernel is $J$ (cf. the second part of Theorem 2.2 (b)). 

Let $z:\Spec(W(k))\to\Spec(S)$ be another morphism and let $w:S\twoheadrightarrow W(k)$ be the $W(k)$-epimorphism that defines $z$. For $i\in\{1,\ldots,n\}$, let $y_i:=\sigma(w(x_i))-w(x_i)^p\in pW(k)$. The filtered Dieudonn\'e module of $(q\circ z)^*(\grD)$ is of the form $(M,F^1_z,\phi)$, where the lift $F^1_z$ of $(M,\phi)$ is computed as follows. Let $u_z:=u_{\text{univ}}\circ z\in U(W(k))\leqslant \pmb{\GL}_M(W(k))$. Let $\phi_z:=u_z\phi$ and let $\tilde\phi_z:({1\over p}F^1+M)\otimes_{W(k)} {}_{\sigma} W(k)\arrowsim M$ be the $W(k)$-linear isomorphism defined by $(M,\phi_z)$. If $\tilde\phi:({1\over p}F^1+M)\otimes_{W(k)} {}_{\sigma} W(k)\arrowsim M$ is the $W(k)$-linear isomorphism defined by $\phi$, then there exists $g_z\in \pmb{\GL}_{({1\over p}F^1+M)\otimes_{W(k)} {}_{\sigma} W(k)}(W(k))$ such that we have $\tilde\phi\circ g_z=\tilde\phi_z$. The {\it correction factor} $g_z$ is computed as follows (see [De, (1.1.3.4)] and [Fa, Sect. 7] for the standard way of pulling back $F$-crystals). For $m\in {1\over p}F^1+M$ we have
$$g_z(m\otimes 1)=\sum_{j_1,j_2,\ldots,j_n\in\dbN\cup\{0\}} (\prod_{i=1}^n \nabla({\partial\over {\partial x_i}})^{j_i})(m\otimes 1)\prod_{i=1}^n {{y_i^{j_i}}\over {j_i!}}.\leqno (14)$$
As $\tilde\phi\circ g_z=\tilde\phi_z$, we have $g_z(F^1\otimes_{W(k)} {}_{\sigma} W(k))=F^1_z\otimes_{W(k)} {}_{\sigma} W(k)$. As $y_1,\ldots,y_n\in pW(k)$, all products ${{y_i^{j_i}}\over {j_i!}}$ belong to $pW(k)$. From this and the fact that $\nabla=\delta_0+\omega$ with $\omega\in E\otimes_{W(k)} \Omega_S^\wedge$, we easily get that (14) implies that there exists an element $e\in pE$ such that $F^1_z=(1_M+e)(F^1)$. Writing $e=e_0+e_{-1}$, where $e_0\in pF^0(E)$ and $e_{-1}\in pN=pF^{-1}(E)$, we have $1_M+e=(1_M+e_{-1})(1_M+e_0)$ and $F^1=(1_M+e_0)(F^1)$. Thus $F^1_z=(1_M+e_{-1})(F^1)$, where $e_{-1}\in pF^{-1}(E)$. Therefore the property 1.6 (iii) holds.

\medskip\noindent
{\bf On 1.6 (iv).} We check that the property 1.6 (iv) holds. Let $D_{W(k)}^z(E)$ be the $p$-divisible group over $W(k)$ that lifts $D$ and that is defined by $z$; its filtered Dieudonn\'e module is $(M,F^1_z,\phi)$. Let $\mu_z:\dbG_m\to \pmb{\GL}_M$ be the inverse of the canonical split  cocharacter of $(M,F^1_z,\phi)$. Let $F^1_z\oplus F^0_z$ be the unique direct sum decomposition normalized by $\im(\mu_z)$. The $W(k)$-modules $\tilde F^0_z(\End(M)):=\Hom(F^0_z,F^1_z)\oplus\End(F^1_z)\oplus \End(F^0_z)$ and $\Hom(F^0_z,F^1_z)$ are the inner conjugates of $\tilde F^0(\End(M)):=\Hom(F^0,F^1)\oplus\End(F^1)\oplus \End(F^0)$ and $\Hom(F^0,F^1)$ (respectively) through $1_M+e_{-1}$. But $E=F^0(E)\oplus F^{-1}(E)$ is fixed under the inner conjugation through $1_M+e_{-1}$. From the last two sentences we get that $F^0(E)=E\cap \tilde F^0(\End(M))=E\cap\tilde F^0_z(\End(M))$ and that $E\cap \Hom(F^0,F^1)=E\cap \Hom(F^0_z,F^1_z)=0$. Thus the following triple $(E,0\subseteq F^0(E)\subseteq E,(0,{1\over p}\phi,\phi))$ is a sub-object of the object 
$$(\End(M),0\subseteq \Hom(F^0_z,F^1_z)\subseteq\tilde F_z^0(\End(M))\subseteq\End(M),(0,{1\over p}\phi,\phi,p\phi))$$ 
of the category of filtered modules $MF_{W,tf}$ over $W(k)$ defined in [Wi, Subsect. 1.1]. From this and the functoriality of canonical split cocharacters (see [Wi, p. 513]), we get that the cocharacter $\mu_{L,z}:\dbG_m\to \pmb{\GL}_{\End(M)}$ defined by $\mu_z$ normalizes $E$. From this and the identity $ E\cap \Hom(F^0_z,F^1_z)=0$, we get that we have a direct sum decomposition $E=F^0(E)\cap F^{-1}_z(E)$, where $F^0(E)=E\cap (\End(F^1_z)\oplus\End(F^0_z))$ and $F^{-1}_z(E):=E\cap \Hom(F^1_z,F^0_z)$. Thus the axiom 1.4 (iii) holds for the pair $(E,D_{W(k)}^z(E))$. As $E=\phi(F^0(E))\oplus p\phi(F^{-1}(E))$, we also have $E=\phi(F^0(E))\oplus p\phi(F^{-1}_z(E))$; therefore the axiom 1.4 (iv) holds for the pair $(E,D_{W(k)}^z(E))$. Thus the property 1.6 (iv) holds.\endproof

\medskip\smallskip\noindent
{\bf 3.6. Proposition (the functorial property).} {\it Let $E$ and $E_1$ be two $W(k)$-submodules of $V_-$ for which the axioms 1.4 (i) and (ii) hold and we have $E\subseteq E_1$. Let $\eps:\dbD^{-1}(E_1)\to \dbD^{-1}(E)$ be the homomorphism defined by the monomorphism $(E,p\phi)\hookrightarrow (E_1,p\phi)$ of Dieudonn\'e modules over $k$. Let $\iota:\dbD^{-1}(E)^{\text{t}}\to \dbD^{-1}(E_1)^{\text{t}}$ be the Cartier dual of $\eps$.

\medskip
{\bf (a)} Then $\scrG_k^E$ is a formal Lie subgroup of $\scrG_k^{E_1}$ and the monomorphism $\scrG_k^E\hookrightarrow\scrG_k^{E_1}$ is isomorphic to the homomorphism of formal Lie groups over $\Spf(k)$ associated to $\iota$ (therefore $\iota$ is a closed embedding homomorphism).

\smallskip
{\bf (b)} Suppose that moreover there exists a $p$-divisible group $D_{W(k)}$ over $W(k)$ that lifts $D$ and such that the axioms 1.4 (iii) and (iv) hold for both pairs $(E,D_{W(k)})$ and $(E_1,D_{W(k)})$. Then $\scrG^E$ is a formal Lie subgroup of $\scrG^{E_1}$ and the monomorphism $\scrG^E\hookrightarrow\scrG^{E_1}$ is isomorphic to the monomorphism of formal Lie groups over $\Spf(W(k))$ associated to the unique monomorphism $\iota_{W(k)}:\dbD^{-1}(E)_{W(k)}^{\text{t}}\hookrightarrow \dbD^{-1}(E_1)_{W(k)}^{\text{t}}$ that lifts $\iota$.}

\medskip
\proof
We can assume that $t\in L_0$ is such that $\phi(t)=t$, $[E,t]=E$, and $[E_1,t]=E_1$ (cf. proof of Lemma 3.1 (c) applied to $(E_1,p\phi)$). Let $n_1:=\dim_k(\nu(E_1))$. Let $S_1:=W(k)[[x_1,\ldots,x_{n_1}]]$. Let $\Phi_{S_1}$, $\Omega_{S_1}^\wedge$, and $d\Phi_{S_1}$ be the analogues of $\Phi_S$, $\Omega_S^\wedge$, and $d\Phi_S$ (respectively) but obtained working with $S_1$ instead of $S$. We identify $\Spec(S)$ with the closed subscheme of $\Spec(S_1)$ defined by the ideal $(x_{n+1},\ldots,x_{n_1})$ of $S_1$. Let $\scrB_1:=(v_1,\ldots,v_{n_1})$ be an $n_1$-tuple of elements of $E_1$ that extends $\scrB=(v_1,\ldots,v_n)$ and such that $\{\nu(v_1),\ldots,\nu(v_{n_1})\}$ is a $k$-basis for $\nu(E_1)$. We consider the deformation package $\dbP(\scrB_1)=(n_1,N_1,n_1,N_1,U_1,U_1,u_{1\text{univ}},\Phi_{N_1},\nabla_1,\grH_{1S/pS},q_{1k},q_{1k},\grK_1,\bar\grK_1)$ of $\scrB_1$. The restriction of $q_1:\Spec(S_1)\to\Spec(R)$ to the closed subscheme $\Spec(S)$ of $\Spec(S_1)$, is the morphism $q:\Spec(S)\to\Spec(R)$; this is so as the $W(k)$-epimorphism $S_1\twoheadrightarrow S$ is compatible with the Frobenius lifts and as $u_{1\text{univ}}\in \pmb{\GL}_M(S_1)$ restricted to $\Spec(S)$ is $u_{\text{univ}}\in \pmb{\GL}_M(S)$. Let $\tilde E_1:=E_1\oplus W(k)t$. Let 
$$0\to (E_1\otimes_{W(k)}  S_1,p\phi\otimes\Phi_{S_1},\delta_0)\to (\tilde E_1\otimes_{W(k)}  S_1,p\tilde\Phi_{N_1},\tilde\nabla_1)\to (S_1t,p\phi\otimes\Phi_{S_1},\delta_0)\to 0\leqno (15)$$
and
$$0\to \dbD^{-1}(E_1)^{\text{t}}\times_k S_1/pS_1\to \grL_{1S_1/pS_1}^{\text{t}}\to (\dbQ_p/\dbZ_p)_{S_1/pS_1}\to 0\leqno (16)$$
be the analogues of (6) and (8) (respectively) but obtained working with $\scrB_1$ instead of $\scrB$. 
The restriction of (15) to $\Spec(S/pS)$ is a short exact sequence 
$$0\to (E_1\otimes_{W(k)}  S,p\phi\otimes\Phi_{S},\delta_0)\to (\tilde E_1\otimes_{W(k)}  S,p\tilde\Phi_{N_1},\tilde\nabla_1)\to (St,p\phi\otimes\Phi_{S},\delta_0)\to 0$$
which is the pull forward of (6) via the monomorphism 
$$(E\otimes_{W(k)}  S,p\phi\otimes\Phi_{S},\delta_0)\hookrightarrow (E_1\otimes_{W(k)}  S,p\phi\otimes\Phi_{S},\delta_0)$$ 
of $F$-crystals over $S/pS$ (this is so as $p\tilde\Phi_{N_1}$ and $p\tilde\Phi_N$ act in the same way on $E_1\otimes_{W(k)}  S$; the connection on $\tilde E_1\otimes_{W(k)}  S$ induced by $\tilde\nabla_1$ is $\tilde\nabla$, cf. Theorem 2.2 (a)). Thus the restriction of (16) to $\Spec(S/pS)$ is a short exact sequence
$$0\to \dbD^{-1}(E_1)^{\text{t}}\times_k S/pS\to \grL_{1S/pS}^{\text{t}}\to (\dbQ_p/\dbZ_p)_{S/pS}\to 0$$
which is the pull forward of (8) via the homomorphism $\tau:\dbD^{-1}(E)^{\text{t}}\to \dbD^{-1}(E_1)^{\text{t}}$. As the addition laws on $\scrG_k^E$ and $\scrG_k^{E_1}$ are defined by adding short exact sequences, we get that the natural inclusion morphism $\scrG_k^E\hookrightarrow\scrG_k^{E_1}$ is a homomorphism of Lie groups; more precisely, it is the homomorphism of formal Lie groups over $\Spf(k)$ associated to $\tau$ (cf. the functorial part of the third paragraph of Section 1). Thus (a) holds. 

The proof of (b) is the same: one only has to use filtered $F$-crystals instead of $F$-crystals (the proof of Lemma 3.1 (a) adapts to show that $\tau_{W(k)}$ exists and it is unique).\endproof
 
\medskip\smallskip\noindent
{\bf 3.7. Corollary.} {\it Let $E$ be a $W(k)$-submodules of $V_-$ for which the axioms 1.4 (i) and (ii) hold. Then there exists a short exact sequence $0\to D_2\to D\to D_1\to 0$ of $p$-divisible groups over $k$ that depends on $E$ and such that $\scrG_k^E$ is naturally a formal Lie subgroup of the formal Lie group $Ext^1(D_1,D_2)$ over $\Spf(k)$.}

\medskip
\proof
We consider a short exact sequence $0\to (M_1,\phi)\to (M,\phi)\to (M_2,\phi)\to 0$ such that $E$ is contained in the direct summand $\Hom(M_2,M_1)=\Hom(M/M_1,M_1)$ of $\End(M)$, cf. Lemma 3.2 (c.ii). Let $0\to D_2\to D\to D_1\to 0$ be the corresponding short exact sequence of $p$-divisible groups over $k$.

Let $D':=D_1\times_k D_2$; its Dieudonn\'e module is $(M',\phi')=(M_1,\phi)\oplus (M_2,\phi)$. We identify $\Hom(M/M_1,M_1)=\Hom(M_2,M_1)=\{x\in \End(M')|x(M_2)=0,\; x(M')\subseteq M_1\}$. Let $E'$ be $E$ but viewed as a $W(k)$-submodule of $\Hom(M_2,M_1)\subseteq \End(M')$. 

Let $\tilde D:=D\times_k D'$; its Dieudonn\'e module is $(\tilde M,\tilde\phi)=(M,\phi)\oplus (M',\phi')$. We can speak about the diagonal embedding $E\cong \tilde E\hookrightarrow E\oplus E'$ inside $\End(\tilde M)$. Let $\scrD_k'$ and  $\tilde\scrD_k'$ be the formal deformation spaces of $D'$ and $\tilde D$ (respectively). We identify naturally $\scrD_k\times_{\Spf(k)} \scrD_k'$ with a formal subspace of $\tilde\scrD_k'$. Let $\scrG_k^E$, $\scrG_k^{'E'}$, and $\tilde \scrG_k^{\tilde E}$ be the formal Lie groups associated to the pairs $(D,E)$, $(D',E')$, and $(\tilde D,\tilde E)$ (respectively). To the diagonal embedding $\tilde E\hookrightarrow E\oplus E'$ inside $\End(\tilde M)$ corresponds (cf. Proposition 3.6 (a)) formal closed embeddings 
$$\tilde \scrG_k^{\tilde E}\hookrightarrow \scrG_k^E\times_{\Spf(k)}  \scrG_k^{'E'}\hookrightarrow\scrD_k\times_{\Spf(k)} \scrD_k'\hookrightarrow\tilde\scrD_k',$$ 
the first one being also a homomorphism. We know that $\scrG_k^{'E'}$ is a formal Lie subgroup of $Ext^1(D_1,D_2)$, cf. Example 3.4. The Corollary follows from the last sentence and the fact that the natural projection homomorphisms $\tilde \scrG_k^{\tilde E}\to \scrG_k^E$ and $\tilde \scrG_k^{\tilde E}\to \scrG_k^{'E'}$ are isomorphisms (one can easily check this at the level of tangent spaces, cf. property 1.5 (ii)).
\endproof

\bigskip
\noindent
{\boldsectionfont 4. Sign $p$-divisible groups} 
\bigskip

We define the positive and the negative $p$-divisible groups of $D$ and we list elementary properties of them (see Subsection 4.1). Their quasi-factors, strings, and slices are introduced in Subsection 4.2. In Subsection 4.3 we combine the Basic Theorem A with Subsection 4.2 to define quasi-factors, strings, and slices of $\scrD_k$ itself. 

For a Newton polygon slope $\alpha\in\scrW$ of $D$ let $r_{\alpha}:=\dim_{B(k)}(W(\alpha))\in\dbN$. Let $\scrX:=\{(\alpha,\beta)\in\scrW\times\scrW|\alpha<\beta\}$. Let $s_1,s_2:\scrX\to\scrW$ be the two projections. We fix a subset $\scrY$ of $\scrX$. We say $\scrY$ is a {\it square zero} subset of $\scrX$, if the intersection $s_1(\scrY)\cap s_2(\scrY)$ is empty. 

\medskip\smallskip\noindent
{\bf 4.1. Lattices.} Let $O_+:=\{x\in V_+|\phi^q(x)\in V_+\;\forall q\in\dbN\}$; it is a nilpotent subalgebra of $V_+$. Let $O_-:=\{x\in V_-|\phi^{-q}(x)\in V_-\;\forall q\in\dbN\}$; it is a nilpotent subalgebra of $V_-$.

\medskip\noindent
{\bf 4.1.1. Lemma.} {\it  The $W(k)$-module $O_-$ is the largest $W(k)$-submodule of $V_-$ such that the pair $(O_-,p\phi)$ is a Dieudonn\'e module over $k$. The $W(k)$-submodule $O_+$ of $V_+$ is the largest $W(k)$-submodule of $V_+$ with the property that the pair $(O_+,\phi)$ is a Dieudonn\'e module.}

\medskip
\proof 
The first part is a particular case of Lemma 3.1 (a). The proof of the second part is similar to the proof of Lemma 3.1 (a) and thus it is left as an exercise.\endproof

\medskip\noindent
{\bf 4.1.2. The trace map.} Let $\Tr:\End(M[{1\over p}])\times\End(M[{1\over p}])\to B(k)$ be the trace map: if $x,y\in\End(M[{1\over p}])$, then $\Tr(x,y)$ is the trace of the endomorphism $xy\in\End(M)$. For $x,y\in \End(M[{1\over p}])$ we have $\Tr(\phi(x),\phi(y))=\sigma(\Tr(x,y))$ i.e., $\Tr$ is invariant under $\phi$. 

Let $\scrL_+=\scrL_+(\scrY):=\oplus_{(\alpha,\beta)\in\scrY}\Hom(W(\alpha),W(\beta))$ and $\scrV_+=\scrV_+(\scrY):=V_+\cap \scrL_+$. Let $\scrL_-=\scrL_-(\scrY):=\oplus_{(\alpha,\beta)\in\scrY}\Hom(W(\beta),W(\alpha))$ and $\scrV_-=\scrV_-(\scrY):=V_-\cap \scrL_-$. Let $\scrO_+=\scrO_+(\scrY):=O_+\cap \scrV_+$ and $\scrO_-=\scrO_-(\scrY):=O_-\cap \scrV_-$. The largest Dieudonn\'e module contained in $(\scrV_+,\phi)$ (resp. in $(\scrV_-,p\phi)$) is $(\scrO_+,\phi)$ (resp. is $(\scrO_-,p\phi)$). 

If $(\alpha_1,\beta_1),(\alpha_2,\beta_2)\in\scrW\times\scrW$ are two pairs such that $(\alpha_1,\beta_1)\neq (\beta_2,\alpha_2)$, then we have $\Tr(\Hom(W(\alpha_1),W(\beta_1)),\Hom(W(\alpha_2),W(\beta_2)))=0$. This implies that:

\medskip
{\it The map $\Tr$ restricts to a perfect bilinear map $\Tr_{\scrY}:\scrL_+\times \scrL_-\to B(k)$ (i.e., the $B(k)$-linear map $\scrL_+\to \Hom(\scrL_-,B(k))$ induced by $\Tr_{\scrY}$ is an isomorphism).}

\medskip
By the {\it dual} of a $W(k)$-lattice $\star$ of $\scrL_+$ (resp. of $\scrL_-$), we mean the unique $W(k)$-lattice $\dbB(\ast)$ of $\scrL_-$ (resp. of $\scrL_+$) such that $\Tr_{\scrY}$ restricts to a perfect map $\star\times \dbB(\ast)\to W(k)$ (resp. $\dbB(\ast)\times \ast\to W(k)$). Let $\scrV_{+,-}=\scrV_{+,-}(\scrY):=\dbB(\scrV_-)\subseteq \scrL_+$ and $\scrV_{-,+}=\scrV_{-,+}(\scrY):=\dbB(\scrV_+)\subseteq \scrL_-$. As $\Tr(\scrV_+,\scrV_-)\subseteq W(k)$, we have inclusions
$$\scrV_+\subseteq \scrV_{+,-}\;\;\text{and}\;\; \scrV_-\subseteq\scrV_{-,+}.\leqno (17)$$
\noindent
{\bf 4.1.3. External mixed sign level modules.} Let $\scrO_{+,-}=\scrO_{+,-}(\scrY)$ be the smallest $W(k)$-submodule of $L_+$ such that $\scrV_{+,-}\subseteq \scrO_{+,-}$ and the pair  $(\scrO_{+,-},\phi)$ is a Dieudonn\'e module over $k$. Let $\scrO_{-,+}=\scrO_{-,+}(\scrY)$ be the smallest $W(k)$-submodule of $L_-$ such that $\scrV_{-,+}\subseteq \scrO_{-,+}$ and the pair $(\scrO_{-,+},p\phi)$ is a Dieudonn\'e module over $k$. We have (cf. (17)):
$$\scrO_+\subseteq \scrO_{+,-}\;\;\text{and}\;\;\scrO_-\subseteq\scrO_{-,+}.\leqno (18)$$
We call $O_{+,-}:=\scrO_{+,-}(\scrX)$ (resp. $O_{-,+}:=\scrO_{-,+}(\scrX)$) the {\it external mixed positive} (resp. {\it negative}) {\it level module} of $(M,\phi)$. One can define the {\it external positive} (resp. {\it negative}) {\it level module} $O_{+,+}$ (resp. $O_{-,-}$) of $(M,\phi)$ and the {\it internal mixed positive} (resp. {\it negative}) {\it level module} of $(M,\phi)$; as they are not used in what follows, we will only say here that $O_{+,+}$ is the smallest $W(k)$-submodule of $L_+$ that contains $V_+$ and such that the pair  $(O_{+,+},\phi)$ is a Dieudonn\'e module over $k$.

\medskip\noindent
{\bf 4.1.4. Definition.} Let $D_+=D_+(\scrY)$ (resp. $D_-=D_-(\scrY)$) be a $p$-divisible group over $k$ whose Dieudonn\'e module is isomorphic to $(\scrO_+,\phi)$ (resp. to $(\scrO_-,p\phi)$). Let $D_{+,-}=D_{+,-}(\scrY)$ (resp. $D_{-,+}=D_{-,+}(\scrY)$) be a $p$-divisible group over $k$ whose Dieudonn\'e module is isomorphic to $(\scrO_{+,-},\phi)$ (resp. to $(\scrO_{-,+},p\phi)$). We call $Z_+:=D_+(\scrX)$ (resp. $Z_-:=D_-(\scrX)$) the {\it positive} (resp. the {\it negative}) $p$-divisible group of either $D$ or $(M,\phi)$. We call $Z_{+,-}:=D_{+,-}(\scrX)$ (resp. $Z_{-,+}:=D_{-,+}(\scrX)$) the {\it external mixed positive} (resp. {\it negative)} {\it $p$-divisible group} of either $D$ or $(M,\phi)$.
 
\medskip\noindent
{\bf 4.1.5. Standard homomorphisms.} To the inclusions of (18) correspond isogenies $D_{+,-}\to D_+$ and $D_{-,+}\to D_-$. As the trace map $\Tr_{\scrY}$ is perfect and as $\Tr$ is invariant under $\phi$, the map $\star\to \dbB(\star)$ induces a bijection between:

\medskip
--  the set of $W(k)$-lattices $\star$ of $\scrL_+$ contained in (resp. containing) $\scrV_+$ with the property  that the pair $(\star,\phi)$ is a Dieudonn\'e module, and

\smallskip
-- the set of $W(k)$-lattices $\ast$ of $\scrL_-$ which contain (resp. are contained in) $\scrV_{-,+}$ with the property that the pair $(\ast,p\phi)$ is a Dieudonn\'e module. 

\medskip\noindent
This bijection is decreasing (i.e., it reverses inclusions). From the last two sentences we get that $\dbB(\scrO_+)=\scrO_{-,+}$ and that $\dbB(\scrO_-)=\scrO_{+,-}$. The perfect bilinear map $(\scrO_{+,-},\phi)\times (\scrO_-,p\phi)\to (W(k),p\sigma)$ defined by $\Tr_{\scrY}$ allows us to identify $D_{+,-}$ with the Cartier dual of $D_-$. Similarly, the perfect bilinear map $(\scrO_+,\phi)\times (\scrO_{-,+},p\phi)\to (W(k),p\sigma)$ defined by $\Tr_{\scrY}$ allows us to identify $D_+$ with the Cartier dual of $D_{-,+}$. 

\medskip\smallskip\noindent
{\bf 4.2. Quasi-factors, strings, and slices.} Let $c_-(\scrY)\in\dbN\cup\{0\}$ be the codimension of $(\scrO_-,p\phi)$. Let $c_-:=c_-(\scrX)$. The pair $(\scrO_-,p\phi)$ is the largest Dieudonn\'e module contained in $(\scrV_-,p\phi)$ and we have $\dim(\nu(\scrO_-))=c_-(\scrY)$, cf. Lemma 3.1 (a) and (b). The inclusion map $\oplus_{(\alpha,\beta)\in\scrX} (\scrO_-(\{(\alpha,\beta)\}),p\phi)\hookrightarrow (O_-,p\phi)$ is an isogeny. Thus we have an identity
$$c_-=\sum_{(\alpha,\beta)\in\scrX} c_-(\{(\alpha,\beta)\})\leqno (19)$$ 
and the inclusion map $(O_{+,-},\phi)\hookrightarrow \oplus_{(\alpha,\beta)\in\scrX} (\scrO_{+,-}(\{(\alpha,\beta)\}),\phi)$ is also an isogeny. The unique Newton polygon slope of $(\Hom(W(\beta),W(\alpha)),p\phi)$ is $1-\beta+\alpha$ and moreover $\dim_{B(k)}(\Hom(W(\beta),W(\alpha))=r_{\alpha}r_{\beta}$. As $\beta-\alpha=1-(1-\beta+\alpha)$, we conclude that
$$c_-(\{(\alpha,\beta)\})=r_{\alpha}r_{\beta}(\beta-\alpha).\leqno (20)$$ 
\indent
We call $D_-(\{(\alpha,\beta)\})$ (resp. $D_{+,-}(\{(\alpha,\beta)\})$) a {\it quasi-factor} of $Z_-$ (resp. of $Z_{+,-}$). For $\alpha\in\scrW$, let $\scrX^{\alpha}:=\{(\beta,\alpha)\in\scrX|\beta\in\scrW\}$ and $\scrX_{\alpha}:=\{(\alpha,\beta)\in\scrX|\beta\in\scrW\}$. We call $D_-(\scrX^{\alpha})$ (resp. $D_-(\scrX_{\alpha})$) the {\it lower} (resp. {\it upper}) {\it $\alpha$-string} of $Z_-$; similarly we call $D_{+,-}(\scrX^{\alpha})$ (resp. $D_{+,-}(\scrX_{\alpha})$) the {\it lower} (resp. {\it upper}) {\it $\alpha$-string} of $Z_{+,-}$. 

We call $D_-=D_-(\scrY)$ (resp. $D_{+,-}=D_{+,-}(\scrY)$) the $\scrY$-slice of $Z_-$ (resp. of $Z_{+,-}$).

\medskip\noindent
{\bf 4.2.1. Fact.} {\it Suppose that $D$ is a direct sum of isoclinic $p$-divisible groups over $k$ (i.e., that $M=\oplus_{\alpha\in\scrW} M\cap W(\alpha)$). Then for each subset $\scrY$ of $\scrX$, we have the following two direct product decompositions $D_-=\prod_{(\alpha,\beta)\in\scrY} D_-(\{(\alpha,\beta)\})$ and $D_{+,-}=\prod_{(\alpha,\beta)\in\scrY} D_{+,-}(\{(\alpha,\beta)\})$.}

\medskip
\proof
We have a direct sum decomposition $V_-=\oplus_{(\alpha,\beta)\in\scrX} \Hom(M\cap W(\beta),M\cap W(\alpha))$ that is stable under $p\phi$. From this we get that $\scrO_-=\oplus_{(\alpha,\beta)\in\scrY} \scrO_-(\{(\alpha,\beta)\})$, cf. definition of $\scrO_-$. Thus $D_-=\prod_{(\alpha,\beta)\in\scrY} D_-(\{(\alpha,\beta)\})$. Taking Cartier duals we get the second product decomposition.\endproof

\medskip
Here is a key particular case of the Basic Theorem A and of the Proposition 3.6 (a).

\medskip\smallskip\noindent
{\bf 4.3. Basic Corollary.} {\it Let $\scrY$ be a square zero subset of $\scrX$. Then  there exists a canonical formal subscheme $\scrG_k(\scrY)$ of $\scrD_k$ that has the following three properties:

\medskip
{\bf (i)} it has a canonical structure of a formal Lie group over $\Spf(k)$ isomorphic to the formal Lie group over $\Spf(k)$ associated to $D_{+,-}=D_{+,-}(\scrY)$;

\smallskip
{\bf (ii)} its tangent space is $\nu(\scrO_-)=\nu(\scrO_{-}(\scrY))$; 

\smallskip
{\bf (iii)} the restriction of $\grD^{\text{f}}$ to $\scrG_k(\scrY)$ is geometrically constant.

\medskip
Let $\scrY_1$ be a subset of $\scrY$. Let $D_{+,-}(\scrY_1)\hookrightarrow D_{+,-}$ be the monomorphism whose Cartier dual is defined by the monomorphism $(\scrO_-(\scrY_1),p\phi)\hookrightarrow (\scrO_-,p\phi)$ of Dieudonn\'e modules over $k$. Then $\scrG_k(\scrY_1)$ is a formal Lie subgroup of $\scrG_k(\scrY)$ and the monomorphism $\scrG_k(\scrY_1)\hookrightarrow\scrG_k(\scrY)$ is isomorphic to the monomorphism of formal Lie groups over $\Spf(k)$ associated to the monomorphism $D_{+,-}(\scrY_1)\hookrightarrow D_{+,-}$.}
  
\medskip
\proof
Let $E:=\scrO_-$. The axiom 1.4 (i) holds, cf. the very definition of $\scrO_-$. If $(\alpha,\beta)\in\scrY$, then as $\scrY$ is a square zero subset of $\scrX$, we have $\scrY\cap\scrX_{\beta}=\emptyset$. Thus $\Hom(W(\beta),W(\alpha))E=0$. This implies that $E^2=0$ i.e., the axiom 1.4 (ii) holds. As the axioms 1.4 (i) and (ii) hold, the Corollary follows from the Basic Theorem A and the Proposition 3.6 (a).\endproof

\medskip\noindent
{\bf 4.3.1. Definition.} We call $\scrG_k(\{(\alpha,\beta)\})$'s (resp. $\scrG_k(\scrY)$'s with $\scrY$ a square zero subset of $\scrX$) the {\it negative quasi-factors} (resp. {\it slices})  of $\scrD_k$. We call $\scrG_k(\scrX^{\alpha})$'s (resp. $\scrG_k(\scrX_{\alpha})$) the {\it negative lower} (resp. {\it upper}) {\it $\alpha$-string} of $\scrD_k$.

\medskip\noindent
{\bf 4.3.2. Example.}  Let $s\in\dbN$ be such that there exists a chain of strict inclusions
$$\scrY_s\subsetneq \scrY_{s-1}\subsetneq\cdots \subsetneq\scrY_1\subsetneq \scrY$$
of non-empty subsets of $\scrX$. If $\scrY$ is a square zero subset of $\scrX$, then we have a chain 
$$\scrG_k(\scrY_s)\hookrightarrow \scrG_k(\scrY_{s-1})\hookrightarrow \cdots \hookrightarrow \scrG_k(\scrY_1)\hookrightarrow  \scrG_k(\scrY)$$
of strict monomorphisms of formal Lie groups over $\Spf(k)$ (cf. Corollary 4.3). 

Let $m>1$ be the number of Newton polygon slopes of $D$ (i.e., of elements of $\scrW$). We write $\scrW=\{\alpha_1,\ldots,\alpha_m\}$ in such a way that $\alpha_1<\cdots<\alpha_m$. If $l\in\{1,\ldots,m-1\}$, then the subset $\scrY[l]=\{(\alpha_i,\alpha_j)|1\le i\le l<j\le m\}$ is a square zero subset and has $l(m-l)$ elements. Note that $\scrG_k(\scrY[1])$ (resp. $\scrG_k(\scrY[m-1])$) is a negative upper $\alpha_1$-string (resp. lower $\alpha_m$-string) of $\scrD_k$ that contains exactly $m-1$ negative quasi-factors of $\scrD_k$. 

Each $\scrY$ has at most $[{m\over 2}](m-[{m\over 2}])$ elements. Therefore we have $s\le [{m\over 2}](m-[{m\over 2}])-1$ and the equality can always take place.

\bigskip
\noindent
{\boldsectionfont 5. Applications to Traverso stratifications}
\bigskip

In Subsections 5.1 to 5.3 we define the most general local Traverso strata and list basic properties of them. Example 5.4 pertains to the principally quasi-polarized context.

Let $\scrA$ be a smooth scheme over $\dbF_p$ which is equidimensional of dimension $cd$. We assume that there exists a $p$-divisible group $\grD_{\scrA}$ over $\scrA$ of codimension $c$ and dimension $d$ which is a versal deformation at all closed points of $\scrA$ with values in finite fields and for which there exists a point $y:\Spec(k)\to\scrA$ such that $D=y^*(\grD_{\scrA})$. This implies that there exists a unique morphism $\varrho:\Spec(R/pR)\to\scrA$ through which $y$ factors and for which we have an identity $\grD\times_R R/pR=\varrho^*(\grD_{\scrA})$ that lifts the identity $D=y^*(\grD_{\scrA})$. The morphism $\varrho$ gives birth to a formally \'etale morphism $\varrho_k:\Spec(R/pR)\to\scrA_k$ of $k$-schemes. Fixing $c$ and $d$, there exist pairs $(\scrA,\grD_{\scrA})$ such that $y$ exists for each $D$ (cf. [NVW, Prop. 2.3]). 

It is known that there exists a unique locally closed, reduced subscheme $\grs_D$ of $\scrA_k$ which has the following property (cf. either [Va2, Thm. 5.3.1 (b)] or [Va3, Thm. 1.2 (d)]): A geometric point $z:\Spec(K)\to\scrA_k$ factors through $\grs_D$ if and only if $z^*(\scrD_{\scrA})$ is isomorphic to $D_K$. The $k$-scheme $\grs_D$ is regular and equidimensional, cf. [Va2, Thm. 5.3.1 (b)] or [Va3, Thm. 1.2 (a)]. We recall that $\grs_D$ is a {\it stratum} of the Traverso stratification of $\scrA_k$ (or of $\scrA$ itself) introduced in [Va2, Thm. 5.3.1 and Subsubsect. 5.3.2] (cf. also [Va3, Subsubsect. 4.2.4]). Let $\scrC_k$ be the formal subscheme of $\scrD_k$ defined by the closed subscheme $\varrho^*(\grs_D)$ of $\Spec(R/pR)$; it is smooth as the morphism $\varrho_k:\Spec(R/pR)\to\scrA_k$ of $k$-schemes is formally \'etale. Let $\tau_D:=\dim(\scrC_k)=\dim(\grs_D)$. 

Let $\scrW$, $\scrX$, $O_-$, $r_{\alpha}$'s, and $c_-$ be as in Section 4. Thus $c_-=\sum_{(\alpha,\beta)\in\scrX} r_{\alpha}r_{\beta}(\beta-\alpha)$, cf. (19) and (20). It is known that 
$$\tau_D=c_-=\sum_{(\alpha,\beta)\in\scrX} r_{\alpha}r_{\beta}(\beta-\alpha).\leqno (21)$$ 
\indent
Formula (21) was first stated in [Tr, Sect. 1, p. 48]: the codimension $cd-\tau_D$ of $\grs_D$ in $\scrA_k$ is $cd-{1\over 2}\sum_{(\alpha,\beta)\in\scrW\times\scrW} |c_{\alpha}d_{\beta}-c_{\beta}d_{\alpha}|=cd-\sum_{(\alpha,\beta)\in\scrX} r_{\alpha}r_{\beta}(\beta-\alpha)$, where $c_{\alpha}:=(1-\alpha)r_{\alpha}$ and $d_{\alpha}:=\alpha r_{\alpha}$ are the codimension and the dimension (respectively) of $(W(\alpha),\phi)$. Strictly speaking, Traverso used a codimension in a certain orbit space context (see [Tr, Thm. 2]) but this is the same as the codimension of $\grs_D$ in $\scrA_k$. We call Formula (21) (or its codimension form) as {\it Traverso's formula}. Recently, the Traverso formula was proved in [Oo] and [Va3, Thm. 1.2 (f)]. The proofs in the last two references are very computational and rely on the fact that $\dim(\grs_D)$ is an isogeny invariant of $D$. 

But the inequality $\tau_D\ge c_-$ can be proved directly using {\it only} local  deformation theory (and thus {\it without} using any global or isogeny arguments). Even more, one can describe $\scrC_k$ directly as follows. Let $(n,S,\Phi_S,\Omega_S,d\Phi_S)$ be as in the beginning of Section 2.

\medskip\smallskip\noindent
{\bf 5.1. Description of $\scrC_k$.} Let $\dbP(\scrB)=(n,N,n,N,U,U,u_{\text{univ}},\Phi_N,\nabla,\grH_{S/pS},q_k,q_k,\grK,\bar\grK)$ be a deformation package with $N\subseteq O_-$ and $n=c_-=\dim_k(\nu(O_-))$ (cf. Lemma 3.1 (b)). 
The $p$-divisible group $\grH_{S/pS}$ over $S/pS$ is geometrically constant, cf. Lemma 3.1 (b) applied to $N\subseteq E=O_-$. Thus the closed subscheme $q_k:\Spec(S/pS)\hookrightarrow\Spec(R/pR)$ is contained in $\varrho^*(\grs_D)$. As $S/pS=k[[x_1,\ldots,x_n]]$, we get $\tau_D=\dim(\varrho^*(\grs_D))\ge n=c_-$. 

Due to the Formula (21), the formal closed embedding $\scrC_k\hookrightarrow \scrD_k=\Spf(R/pR)$ gets identified with the formal closed embedding $q_k^{\text{f}}:\Spf(S/pS)\hookrightarrow\Spf(R/pR)$.

\medskip\smallskip\noindent
{\bf 5.2. Local Traverso strata.} Let $G$ be a flat, closed subgroup scheme of $\pmb{\GL}_M$ such that the following two properties hold:

\medskip
{\bf (i)} its generic fibre $G_{B(k)}$ is connected and has a Lie algebra $\Lie(G_{B(k)})$ that is left invariant under $\phi$;

\smallskip
{\bf (ii)} there exists a cocharacter $\mu:\dbG_m\to\pmb{\GL}_M$ that normalizes $G$ and that gives birth to a direct sum decomposition $M=F^1\oplus F^0$ such that $F^1/pF^1=\bar F^1$ and $\mu$ acts on each $F^i$ as the $-i$-th power of the identity character of $\dbG_m$.

\medskip
Based on the functorial aspects of canonical split cocharacters (see [Wi, p. 413]), we can assume that $\mu:\dbG_m\to\pmb{\GL}_M$ is the canonical split cocharacter of the filtered Dieudonn\'e module $(M,F^1,\phi)$. The largest $W(k)$-submodule of $\End(M)\cap\Lie(G_{B(k)})$ on which $\dbG_m$ acts through $\mu_L:\dbG_m\to \pmb{\GL}_{\End(M)}$ via the identity character of $\dbG_m$, is 
$$N_G(\mu):=\Hom(M/F^0,F^0)\cap\Lie(G_{B(k)})$$ 
(the intersection being taken inside $\End(M[{1\over p}])$). 
Let $n_G$ be the rank of $N_G(\mu)$. We have $n_G=\dim_k(\nu(N_G(\mu)))$.

Let $V_-(G):=V_-\cap\Lie(G_{B(k)})$. Let $O_-(G):=O_-\cap V_-(G)$. Let $D_-(G)$ be the $p$-divisible group over $k$ whose Dieudonn\'e module is $(O_-(G),p\phi)$. Let $c_-(G)\in\dbN\cup\{0\}$ be the codimension of $D_-(G)$; thus $c_-(G)=\dim_k(\nu(O_-(G)))$ (cf. Lemma 3.1 (b)). 

Let $S_G:=W(k)[[x_1,\ldots,x_{n_G}]]$. Let $\Phi_{S_G}$ be the Frobenius lift of $S_G$ that is compatible with $\sigma$ and that takes $x_i$ into $x_i^p$ for all $i\in\{1,\ldots,n_G\}$. Let $\scrB_G=(v_1,\ldots,v_{n_G})$ be an $n_G$-tuple of elements of $\End(M)\cap\Lie(G_{B(k)})$ such that $\{\nu(v_1),\ldots,\nu(v_{n_G})\}$ is a $k$-basis for $\nu(N_G(\mu))$. Let $N_G:=\oplus_{i=1}^{n_G} W(k)v_i$; it is a direct summand of $\End(M)$. The rule $U_G(A):=\{x\in 1_{M\otimes_{W(k)} A}+N_G\otimes_{W(k)} A|\text{det}(x)\,\,\text{is}\,\,\text{an}\,\,\text{invertible}\,\,\text{element}\,\,\text{of}\,\,A\}$ defines a closed subscheme $U_G$ of $\pmb{\GL}_M$. 

We will choose $\scrB_G$ such that $U_G$ is a closed subscheme of $G$. If $N_G^2=0$, then $U_G$ is a closed subgroup scheme of $\pmb{\GL}_M$ whose Lie algebra is exactly $N_G$; thus the generic fibre of $U_G$ is a closed subgroup of $G_{B(k)}$ (cf. [Bo, Ch. II, Subsect. 7.1]) and therefore $U_G$ is a closed subscheme of $G$. For instance, we can take $\scrB_G$ such that $N_G=N_G(\mu)$. 

Let $\dbP(\scrB_G)=(n_G,N_G,n_G,N_G,U_G,U_G,u_{G\text{univ}},\Phi_{N_G},\nabla_G,\grH_{S_G/pS_G},q_{Gk},q_{Gk},\grK_G,\bar\grK_G)$ be the deformation package of $\scrB_G$ obtained working with $S_G$ instead of $S$.

Let $\scrD_k(G,\scrB_G)$ be the smooth formal subscheme of $\scrD_k$ defined by the closed embedding $q_{Gk}:\Spec(S_G/pS_G)\hookrightarrow\Spec(R/pR)$; its tangent space is $\nu(N_G)=\nu(N_G(\mu))$ (cf. Lemma 2.3). Let 
$$\scrC_k(G,\scrB_G):=\scrC_k\cap\scrD(G,\scrB_G);$$ 
its tangent space $\grt$ is $\nu(N_G(\mu))\cap \nu(O_-)$. As $\nu(O_-(G))$ is contained in $\nu(N_G(\mu))\cap \nu(O_-)=\nu(\Lie(G_{B(k)})\cap\End(M))\cap \nu(O_-)$, we get:

\medskip
{\it The tangent space $\grt$ of  $\scrC_k(G,\scrB_G)$ has dimension at least equal to $c_-(G)$.}

\medskip
One can check that the following two properties hold:

\medskip
{\bf $I_G$.} The formal scheme $\scrD_k(G,\scrB_G)$ (i.e., $q_{Gk}$) depends only on $(M,\phi,G)$ and not on the choice of $\scrB_G$ (one can check this in a way similar to Proposition 2.6 (a) and Corollary 2.7, to be compared with [Fa, Sect. 7, Rm. iii), p. 136]). 

\smallskip
{\bf $II_G$.} If $G$ is smooth over $W(k)$, then the reduced formal scheme $\scrC_k(G,\scrB_G)_{\text{red}}$ of $\scrC_k(G,\scrB_G)$ is formally smooth (the proof of this is the same as of [Va3, Cor. 4.3 (a)]).

\medskip
To be short, we will not check these two properties here and thus in what follows we will keep using the notations $\scrD_k(G,\scrB_G)$  and $\scrC_k(G,\scrB_G)$. 
Based on $I_G$ and $II_G$, one expects that, under mild assumptions, the formal scheme $\scrC_k(G,\scrB_G)$ is isomorphic to $\Spf(k[[x_1,\ldots,x_{c_-(G)}]])$. This is equivalent to the following two statements:

\medskip\noindent
{\bf (a)} We have $\dim_k(\grt)=c_-(G)$.

\smallskip\noindent
{\bf (b)} There exists a formal closed subscheme of $\scrC_k(G,\scrB_G)$ isomorphic to $\Spf(k[[x_1,\ldots,x_{c_-(G)}]])$.

\medskip
For the remaining part of this Subsection we will concentrate on the statement (a). Subsection 5.3 will be centered on the statement (b). 

\medskip\noindent
{\bf 5.2.1. Definitions.} We call $\scrD_k(G,\scrB_G)$ the {\it formal deformation space} of the triple $(M,\phi,G)$ (with respect to $\scrB_G$). We call $\scrC_k(G,\scrB_G)_{\text{red}}$ (resp. $\scrC_k(G,\scrB_G)$) the {\it local} (resp. {\it potentially non-reduced local}) {\it Traverso stratum} of $(M,\phi,G)$ (with respect to $\scrB_G$).

\medskip\noindent
{\bf 5.2.2. Fact.} {\it The following two properties hold:

\medskip
{\bf (a)} We have $\dim_k(\grt)=c_-(G)$ if and only if $\nu(O_-(G))=\nu(V_-(G))\cap \nu(O_-)$. 

\smallskip
{\bf (b)} If we have a direct sum  decomposition $(V_-,p\phi)=(V_-(G),p\phi)\oplus (V_-(G)^\perp,p\phi)$ of $F$-crystals over $k$, then we have $\nu(O_-(G))=\nu(V_-(G))\cap \nu(O_-)$.}

\medskip
\proof
Part (a) is obvious.  To check (b), let $(O_-(G)^\perp,p\phi)$ be the largest Dieudonn\'e module contained in $(V_-(G)^\perp,p\phi)$. We have a direct sum decomposition $O_-=O_-(G)\oplus O_-(G)^\perp$, cf. the very definition of $O_-$. This implies that $\nu(V_-(G))\cap \nu(O_-)=\nu(O_-(G))$.\endproof

\medskip\noindent
{\bf 5.2.3. Example.} Suppose that there exists a reductive, closed subgroup scheme $H$ of $G$ such that $\Lie(H_{B(k)})$ is stable under $\phi$, we have $V_-(G)\subseteq\Lie(H)$, and the restriction of the trace map $\Tr$ (introduced in Subsubsection 4.1.2) to $\Lie(H)\times \Lie(H)$ is a perfect bilinear form. Let $\Lie(H)^\perp$ be the perpendicular on $\Lie(H)$ with respect to $\Tr:\End(M)\times\End(M)\to W(k)$. Let $V_-(G)^\perp:=V_-\cap \Lie(H)^\perp$. As $V_-(G)=V_-\cap\Lie(H)$ and the direct sum decomposition $\End(M)=\Lie(H)\oplus\Lie(H)^\perp$ of $W(k)$-modules is stable under $p\phi$, we also have a direct sum decomposition $(V_-,p\phi)=(V_-(G),p\phi)\oplus (V_-(G)^\perp,p\phi)$. Thus we have $\dim_k(\grt)=c_-(G)$, cf. Fact 5.2.2.

\medskip\noindent
{\bf 5.2.4. Example.} Suppose $\End(M)\cap\Lie(G_{B(k)})$ is a subalgebra of $\End(M)$. In this case $U_G$ is a closed subscheme of $G$ for each $\scrB_G$. Thus we can assume that $\scrB_1:=(v_1,\ldots,v_{c_-(G)})$ belongs to $O_-(G)^{\times c_-(G)}$.  We consider the deformation package $\dbP(\scrB_1)=(c_-(G),N_1,c_-(G),N_1,U_1,U_1,u_{1\text{univ}},\Phi_{N_1},\nabla_1,\grH_{1S_1/pS_1},q_{1k},q_{1k},\grK_1,\bar\grK_1)$ of $\scrB_1$. The closed embedding $q_{1k}:\Spec(S_1/pS_1)\hookrightarrow\Spec(R/pR)$ is the restriction of $q_{Gk}$ to the closed subscheme $\Spec(S_1/pS_1)$ of $\Spec(S_G/pS_G)$ defined by the ideal $(x_{c_-(G)+1},\ldots,x_{n_G})$ of $S_G/pS_G$. Let $\scrC_k(G,\scrB_1)$ be the smooth formal subscheme of $\scrD_k(G,\scrB_G)$ defined by the closed embedding $q_{1k}:\Spec(S_1/pS_1)\hookrightarrow\Spec(R/pR)$; its dimension is $c_-(G)$. As $\grH_{1S_1/pS_1}$ is geometrically constant (cf. Lemma 3.1 (b)), $\scrC_k(G,\scrB_G)$ contains the closed subscheme $\scrC_k(G,\scrB_1)$. Moreover the following two properties hold:

\medskip
{\bf (a)} If $\dim_k(\grt)=c_-(G)$, then in fact we have $\scrC_k(G,\scrB_G)=\scrC_k(G,\scrB_1)=\Spf(S_1/pS_1)$.

\smallskip
{\bf (b)} If $O_-(G)^2=0$, then Corollary 2.7 implies that $\scrC_k(G,\scrB_1)$ depends only on $O_-(G)$ (i.e., on $(M,\phi,G)$) and the Basic Theorem A implies that $\scrC_k(G,\scrB_1)$ has a canonical formal Lie group structure isomorphic to the formal Lie group over $\Spf(k)$ associated to $D_-(G)^{\text{t}}$.

\medskip\smallskip\noindent
{\bf 5.3. Closed subschemes of $\scrD_k(G,\scrB_G)$ and $\scrC_k(G,\scrB_G)$.} We begin by introducing some closed subgroup schemes of $\pmb{\GL}_M$ and $G$ that can be used to construct explicitly closed subschemes of $\scrD_k(G,\scrB_G)$ and $\scrC_k(G,\scrB_G)$. Let $m>1$ and $\scrW=\{\alpha_1,\ldots,\alpha_m\}$ be as in Example 4.3.2. We have an increasing Newton polygon slope filtration 
$$0\subset M\cap W(\alpha_1)\subset  M\cap [W(\alpha_1)\oplus W(\alpha_2)]\subset\cdots\subset M\cap [\oplus_{i=1}^{m-1}W(\alpha_i)]\subset M.$$
Let $\scrH$ be the parabolic subgroup scheme of $\pmb{\GL}_M$ that is the normalizer of this filtration of $M$. Let $\scrU$ be the unipotent radical of $\scrH$; its Lie algebra is $V_-$. Let $\scrU_{G,B(k)}:=G_{B(k)}\cap \scrU_{B(k)}$; it is a unipotent subgroup of $G_{B(k)}$. Let  $\scrU_G$ be the schematic closure of $\scrU_{G,B(k)}$ in $G$. If $\scrU_G$ is smooth, then we have $\Lie(\scrU_G)=V_-(G)$.

\medskip\smallskip\noindent
{\bf 5.3.1. Lemma.} {\it We assume that either $G$ is a reductive group scheme over $W(k)$ or we have $m\le p$. Then  $\scrU_G$ is a smooth, unipotent, closed subgroup scheme of $G$.}

\medskip
\proof
We first assume that $G$ is a reductive group scheme over $W(k)$. The conditions 5.2 (i) and (ii) imply that the Newton quasi-cocharacter of $\pmb{\GL}_{M[{1\over p}]}$ associated to the direct sum decomposition $M[{1\over p}]=\oplus_{i=1}^m W(\alpha_i)$ normalizes $G_{B(k)}$. From this and the fact that $G$ is a reductive group scheme over $W(k)$, one easily gets that $G\cap\scrH$ is a parabolic subgroup scheme of $G$. The unipotent radical of $G\cap\scrH$ has the same generic fibre as $\scrU_G$ and thus it is equal to $\scrU_G$. Therefore $\scrU_G$ is a smooth, unipotent, closed subgroup scheme of $G$.

If $m\le p$, then $V_-^p=0$ and thus one can use exponential and logarithmic maps to identify the vector group schemes over $W(k)$ associated to $V_-$ and $V_-(G)$ with $\scrU$ and $\scrU_G$ (respectively). Therefore $\scrU_G$ is a smooth, unipotent, closed subgroup scheme of $G$.\endproof

\medskip\smallskip\noindent
{\bf 5.3.2. A more general definition of $\scrD_k(G,\scrB_G)$.} Let $z:\Spec(W(k))\hookrightarrow \Spec(S_G)$ be the closed embedding defined by the ideal $J_G=(x_1,\ldots,x_{n_G})$ of $S_G$. Suppose that there exists a morphism $b_G:\Spec(S_G)\to G$ such that the following two conditions hold:

\medskip
{\bf (i)} the composite $b_G\circ z$ defines the identity section of $G$;

\smallskip
{\bf (ii)} the tangent map of $b_G$ at $z$ is a $W(k)$-linear monomorphism that defines a direct summand $N_{b_G}$ of $\End(M)\cap\Lie(G_{B(k)})$ of rank $n_G$ and with $\dim_k(\nu(N_{b_G}))=n_G$.

\medskip
As in Theorem 2.2, Faltings deformation theory implies that there exists a $p$-divisible group $\grH^{\text{f}}_G$ over $\Spf(S_G)$ whose filtered $F$-crystal over $S_G/pS_G$ is the quadruple 
$$(M\otimes_{W(k)} S_G,F^1\otimes_{W(k)} S_G,b_G(\phi\otimes\Phi_{S_G}),\nabla_G),$$
where the connection $\nabla_G$ is uniquely determined by $b_G(\phi\otimes\Phi_{S_G})$. To it corresponds a formal closed embedding $\tilde q_G^{\text{f}}:\Spf(S_G)\hookrightarrow\Spf(R)$ such that $\grH_G^{\text{f}}=\tilde q_G^{\text{f}*}(\grH^{\text{f}})$ (cf. (ii)) which defines a formal closed subscheme $\tilde q_{Gk}^{\text{f}}:\scrC_k(G,b_G)\hookrightarrow \scrD_k$.  

One can check that the following property holds:

\medskip
{\bf $III_G$.} The formal scheme $\scrD_k(G,b_G)$ (i.e., $\tilde q_{Gk}^{\text{f}}$) depends only on $(M,\phi,G)$ and not on the choice of $b_G$ and is equal to each $\scrC_k(G,\scrB_G)$ (one can check this in a way similar to Proposition 2.6 (a) and Corollary 2.7, to be compared with [Fa, Sect. 7, Rm. iii), p. 136]). 

\medskip
Let $n_1=c_-(G)$ and let $z_1:\Spec(S_1)\hookrightarrow \Spec(S_G)$ be the closed embedding defined by the ideal $(x_{n_1+1},\ldots,x_{n_G})$ of $S_G$. We consider the following two conditions:

\medskip
{\bf (iii)} the composite $b_G\circ z_1:\Spec(S_1)\to G$ factors through $\scrU_G$ and thus through $\scrU$;

\smallskip
{\bf (iv)} condition (iii) holds and the resulting $S_1$-valued point of $\scrU$ is an element of the subgroup $1_{M\otimes_{W(k)} S_1}+O_-\otimes_{W(k)} S_1$ of $\scrU(S_1)$.

\medskip
If (iv) holds, then the restriction of $\grH^{\text{f}}_G$ to $\Spf(S_1/pS_1)$ is geometrically constant (the proof of this is the same as Fact 2.3 and Lemma 3.1 (b)); thus the composite $\Spf(S_1/pS_1)\hookrightarrow\Spf(S_G/pS_G)\hookrightarrow \Spf(R/pR)$ associated to $z_1$ and $\tilde q_{Gk}$ defines a formal closed subscheme $\scrC_k(G,b_G)$ of the intersection $\scrD_k(G,b_G)\cap\scrC_k$. We have the following simple Fact.

\medskip\smallskip\noindent
{\bf 5.3.3. Fact.} {\it If $G$ and $\scrU_G$ are smooth over $W(k)$ (for instance, this holds if $G$ is a reductive group scheme over $W(k)$, cf. Lemma 5.3.1), then there exists a morphism $b_G:\Spec(S_G)\to G$ such that the conditions 5.3.2 (i) to (iii) hold. If moreover, via the exponential map we have $\text{exp}(O_-(G))\subseteq \scrU_G(W(k))$ (for instance, this holds if $m\le p$), then there exists a morphism $b_G:\Spec(S_G)\to G$ such that the conditions 5.3.2 (i) to (iv) hold.}

\medskip
\proof
The first part is trivial. To check the second part, let $\scrU_G^{-}$ be the smooth vector group scheme over $W(k)$ whose Lie algebra is $O_-(G)$. The exponetial map defines a homomorphism $e_G:\scrU_G^{-}\to \scrU_G$ whose generic fibre is an isomorphism. The $W(k)$-linear map $de_G$ is the inclusion $O_-(G)=\Lie(\scrU_G^{-})\hookrightarrow V_-(G)=\Lie(\scrU_G)$. We take $b_G$ such that $b_G\circ z_1$ factors through $e_G$. Property (iv) holds as the composite homomorphism $\scrU_G^{-}(S_1)\to \scrU_G(S_1)\to\scrU(S_1)$ induced by $e_G(S_1)$, factors through $1_{M\otimes_{W(k)} S_1}+O_-\otimes_{W(k)} S_1$.\endproof 
 
\medskip\smallskip\noindent
{\bf 5.4. Example (the principally quasi-polarized context).} Suppose that $c=d$ and that $D$ has a principal quasi-polarization $\lambda$. Let $\psi:M\times M\to W(k)$ be the perfect, alternating form defined by $\lambda$. Suppose that $G=\pmb{\text{Sp}}(M,\psi)$. We have $V_-(G)\subseteq\Lie(G)$ and $\Lie(G_{B(k)})$ is stable under $\phi$. If $F^1$ and $F^0$ are anisotropic with respect to $\psi_M$ and if $F^1/pF^1=\bar F^1$, then $\mu$ factors through $\pmb{\text{GSp}}(M,\psi)$ and thus it normalizes $G$. Either $\scrD_k(G,\scrB_G)$ or $\scrD_k(G,b_G)$ is the smooth formal subscheme of $\scrD_k$ which is the formal deformation space of $(D,\lambda)$ over $\Spf(k)$. Thus $I_G$ and $III_G$ hold. Let $\scrD_k(G)=\scrD_k(G,\scrB_G)=\scrD_k(G,b_G)$; its dimension is ${{d(d+1)}\over 2}$. We recall that $n_1=c_-(G)$. 

In this paragraph we assume $p>2$. Let $(v_1,\ldots,v_{n_1})\in O_-(G)^{n_1}$ be such that $\{\nu(v_1),\ldots,\nu(v_{n_1})\}$ is a $k$-basis of $\nu(O_-(G))$. The Cayley element $\varpi:=(1_{M\otimes_{W(k)} S_1}-\sum_{i=1}^{n_1} v_i\otimes x_i)(1_{M\otimes_{W(k)} S_1}+\sum_{i=1}^{n_1} v_i\otimes x_i)^{-1}$ belongs to $G(S_1)\cap [1_{M\otimes_{W(k)} S_1}+O_-\otimes_{W(k)} S_1]$ and thus also to $\scrU_G(S_1)=G(S_1)\cap \scrU(S_1)$. This implies that we can choose $b_G$ such that $b_G\circ z_1:\Spec(S_1)\to G$ is $\varpi$ and therefore such that the conditions 5.3.2 (i) to (iv) hold.

If $p=2$, then in this paragraph we assume that $m\le 2$ (or that $b_G$ is such that the condition 5.3.2 (iv) holds). From the Fact 5.3.3 and the previous paragraph, we get that regardless of what $p$ is, we can choose $b_G$ such that the conditions 5.3.2 (i) to (iv) hold. Thus from the end of Subsection 5.3.2 we get that $\scrC_k(G)$ contains a formal closed subscheme isomorphic to $\Spf(k[[x_1,\ldots,x_{n_1}]])$. If $p>2$, then the restriction of the trace map $\Tr$ to $\Lie(G)\times \Lie(G)$ is a perfect bilinear form. Thus if $p>2$, then we have $\dim_k(\grt)=n_1$ (cf. Example 5.2.3) and therefore $\scrC_k(G)$ is isomorphic to $\Spf(k[[x_1,\ldots,x_{n_1}]])$. 

We use the notations of the beginning paragraphs of Section 4. Manin's symmetric condition says: for $i\in\{1,\ldots,m\}$ we have $\alpha_i+\alpha_{m+1-i}=1$ and $r_{\alpha_i}=r_{\alpha_{m+1-i}}$. We have 
$$V_-(G)[{1\over p}]=\oplus_{(\alpha,\beta)\in\scrX}  V_-(G)[{1\over p}]\cap [\Hom(W(\beta),W(\alpha))\oplus \Hom(W(1-\alpha),W(1-\beta))].$$ 
If $\alpha+\beta\neq 1$ and $(\alpha,\beta)\in\scrX$, then $\psi_M(W(\beta),W(\alpha))=\psi_M(W(1-\alpha),W(1-\beta))=0$ and it is easy to see that $V_-(G)[{1\over p}]\cap [\Hom(W(\beta),W(\alpha))\oplus \Hom(W(1-\alpha),W(1-\beta))]$ has the same dimension $r_{\alpha}r_{\beta}$ as $\Hom(W(\beta),W(\alpha))$. If $\alpha+\beta=1$ and $(\alpha,\beta)\in\scrX$, then $(W(\beta)\oplus W(\alpha),\psi_M)$ is an alternating space over $B(k)$; thus  
$(V_-(G)[{1\over p}]\cap \Hom(W(\beta),W(\alpha)),p\phi)$ is an isoclinic $F$-isocrystal of Newton polygon slope $1-\beta+\alpha=2\alpha$, rank ${{r_{\alpha}(r_{\alpha}+1)}\over 2}$, and codimension ${{r_{\alpha}(r_{\alpha}+1)}\over 2}(1-2\alpha)$. From the last two sentences and (20), similar to (21) we get:
$$n_1={1\over 2}\sum_{(\alpha,\beta)\in\scrX,\alpha+\beta\neq 1} c_-(\{\alpha,\beta\})+\sum_{(\alpha,1-\alpha)\in\scrX} r_{\alpha}(r_{\alpha}+1)({1\over 2}-\alpha)={1\over 2}c_-+\sum_{(\alpha,1-\alpha)\in\scrX} r_{\alpha}({1\over 2}-\alpha).\leqno (22)$$
Formulas equivalent to (22) were first obtained in [Oo].

\bigskip\noindent
{\bf Acknowledgments.} We would like to thank University of Arizona and Binghamton University for good working conditions. We would like to thank C.-L. Chai for mentioning to us that it would be of interest to simplify [Ch2]. We would also like to thank M.-H. Nicole for few comments and the referee for many valuable suggestions and comments. This project was partially supported by the NSF grant DMS \#0900967.

\bigskip
\references{37}
{\nspace{

\Ref[Bo]
A. Borel,
\sl Linear algebraic groups. Second edition,
\rm Grad. Texts in Math., Vol. {\bf 126}, Springer-Verlag, New York, 1991.

\Ref[BM]  P. Berthelot and W. Messing, 
\sl Th\'eorie de Dieudonn\'e crystalline III, 
\rm in the Grothendieck Festschrift I, Progr. Math., Vol. {\bf 86},  173--247, Birkh\"auser, Basel, 1990.

\Ref[Ch1]
C.-L. Chai,
\sl Hecke orbits on Siegel modular varieties,
\rm  Geometric methods in algebra and number theory,  71--107, Progr. Math., Vol. {\bf 235}, Birkh\"auser Boston, Boston, MA, 2005.

\Ref[Ch2]
C.-L. Chai,
\sl Canonical coordinates on leaves of p-divisible groups: The two-slope case,
\rm 77 pages manuscript dated 01/10/2005.

\Ref[De] 
P. Deligne, 
\sl Cristaux ordinaires et coordon\'ees canoniques, 
\rm Lecture Notes in Math., Vol. {\bf 868},  80--137, Springer-Verlag, 1981.

\Ref[Di] 
J. Dieudonn\'e, 
\sl Groupes de Lie et hyperalg\`ebres de Lie sur un corps de caract\'erisque $p>0$ (VII), 
\rm Math. Annalen {\bf 134} (1957),  114--133.

\Ref[Fa] 
G. Faltings, 
\sl Integral crystalline cohomology over very ramified valuation rings, 
\rm J. Amer. Math. Soc. {\bf 12} (1999), no. 1,  117--144.

\Ref[Fo]
J.-M. Fontaine, 
\sl Groupes $p$-divisibles sur les corps locaux, 
\rm J. Ast\'erisque {\bf 47/48}, Soc. Math. de France, Paris, 1977.

\Ref[HT]
M. Harris and R. Taylor, 
\sl The geometry and cohomology of some simple Shimura varieties,
\rm Annals of Mathematics Studies, Vol. {\bf 151}. Princeton University Press, Princeton, NJ, 2001.

\Ref[Il]
L. Illusie, 
\sl D\'eformations des groupes de Barsotti--Tate (d'apr\`es A. Grothendieck), 
\rm Seminar on arithmetic bundles: the Mordell conjecture (Paris, 1983/84),  151--198, J. Ast\'erisque {\bf 127}, Soc. Math. de France, Paris, 1985.

\Ref[Ka1] 
N. Katz, 
\sl Slope filtration of $F$-crystals, 
\rm Journ\'ees de G\'eom\'etrie alg\'ebrique de Rennes 1978, J. Ast\'erisque, Vol. {\bf 63}, Part 1,  113--163, Soc. Math. de France, Paris, 1979.

\Ref[Ka2] 
N. Katz, 
\sl Appendix to Exp. V, 
\rm Lecture Notes in Math., Vol. {\bf 868},  128--137, Springer-Verlag, 1981.

\Ref[Ma]
J. I. Manin, 
\sl The theory of formal commutative groups in finite characteristic, 
\rm Russian Math. Surv. {\bf 18} (1963), no. 6,  1--83.

\Ref[Me] 
W. Messing, 
\sl The crystals associated to Barsotti--Tate groups, with applications to abelian schemes, 
\rm Lecture Notes in Math., Vol. {\bf 264}, Springer-Verlag, Berlin-New York, 1972.

\Ref[Mi1]
J. S. Milne,
\sl The points on a Shimura variety modulo a prime of good
reduction,
\rm The Zeta functions of Picard modular surfaces,  153--255, Univ. Montr\'eal Press, Montreal, Quebec, 1992.

\Ref[Mi2]
J. S. Milne,
\sl Shimura varieties and motives,
\rm Motives (Seattle, WA, 1991),  447--523, Proc. Sympos. Pure Math., Vol. {\bf 55}, Part 2, Amer. Math. Soc., Providence, RI, 1994.

\Ref[Mo]
B. Moonen,
\sl Linearity properties of Shimura varieties. II,
\rm Compos. Math. {\bf 114}  (1998),  no. 1,  3--35. 

\Ref[NVW]
M.-H. Nicole, A. Vasiu, and T. Wedhorn,
\sl Purity of level $m$ stratifications,
\rm to appear in Ann. Sci. \'Ecole Norm. Sup. (4) {\bf 43} (2010), no. 6. 

\Ref[Oo]
F. Oort,
\sl Foliations in moduli spaces of abelian varieties and dimension of leaves, 
\rm to appear in Algebra, Arithmetic, and Geometry, in honor of Yu. I. Manin, Progr. Math., Birkh\"auser, Basel, 2010.

\Ref[Ta]
J. Tate,
\sl $p$-divisible groups, 
\rm Proc. Conf. Local Fields (Driebergen, 1966), 158--183 Springer, Berlin, 1967.

\Ref[Tr]
C. Traverso,
\sl $p$-divisible groups over fields,
\rm  Symposia Mathematica, Vol. {\bf XI} (Convegno di Algebra Commutativa, INDAM, Rome, 1971),  45--65, Academic Press, London, 1973.

\Ref[Va1]
A. Vasiu, 
\sl Integral canonical models of Shimura varieties of preabelian
 type, 
\rm Asian J. Math. {\bf 3} (1999), no. 2,  401--518.

\Ref[Va2]
A. Vasiu,
\sl Crystalline Boundedness Principle,
\rm Ann. Sci. \'Ecole Norm. Sup. (4) {\bf 39} (2006), no. 2, 245--300.

\Ref[Va3]
A. Vasiu,
\sl Level $m$ stratifications of versal deformations of $p$-divisible groups,
\rm J. Alg. Geom. {\bf 17} (2008), no. 4, 599--641.

\Ref[Va4]
A. Vasiu,
\sl Generalized Serre--Tate ordinary theory,
\rm manuscript math.NT/0208216.

\Ref[Wi] 
J.-P. Wintenberger, 
\sl Un scindage de la filtration de Hodge pour certaines vari\'et\'es alg\'ebriques sur les corps locaux, 
\rm Ann. of Math. (2) {\bf 119}, (1984), no. 3,  511--548.

}}
\noindent
\bigskip
\hbox{Adrian Vasiu,}
\hbox{Binghamton University,}
\hbox{Department of Mathematical Sciences,}
\hbox{P.O. Box 6000}
\hbox{Binghamton, New York 13902-6000, U.S.A.}
\hbox{e-mail: adrian\@math.binghamton.edu}

\enddocument